\documentclass[leqno,12pt]{amsart}
\usepackage{amssymb} 
\theoremstyle{plain} 
\newtheorem{theorem}{\indent\sc Theorem}[section] 

\newtheorem{proposition}[theorem]{\indent\sc Proposition}

\theoremstyle{definition} 
\newtheorem{definition}[theorem]{\indent\sc Definition}
\newtheorem{remark}[theorem]{\indent\sc Remark}

\newtheorem{problem}[theorem]{\indent\sc Problem}
\begin{document}

\title{Painlev\'e scheme \\}
\author{Yusuke Sasano }

\renewcommand{\thefootnote}{\fnsymbol{footnote}}
\footnote[0]{2000\textit{ Mathematics Subjet Classification}.
34M55; 34M45; 58F05; 32S65.}


\begin{abstract}
In this note, we review the notion of Painlev\'e scheme of the sixth Painlev\'e equation from the viewpoint of accessible singular point and its local index in the Hirzebruch surface of degree two ${\Sigma_2}$. The key method is Painlev\'e $\alpha$-method for each accessible singular point. Giving a Painlev\'e scheme in the differential system satisfying certain conditions, we can recover the Painlev\'e VI system with the polynomial Hamiltonian. We also consider the case of the Painlev\'e V,IV and III systems, respectively. Finally, we study non-linear ordinary differential systems in dimension two with only simple accessible singular $(n+2)$-points in the Hirzebruch surface of degree $n$; ${\Sigma_n}$. This equation has symmetry of symmetric group of degree $n+2$.
\end{abstract}
\maketitle

\section{Introduction}

For a linear differential equation of Fuchs type, we can make a Riemann scheme. This is the pair of singularity and local exponent. Conversely, by giving the Riemann scheme satisfying the Fuchs relation, we can recover a linear differential equation of Fuchs type.

Painlev\'e equations are the second-order non-linear ordinary differential equations. Here, we consider the following problem.

\begin{problem}
Can we construct a generalization of the Riemann scheme for each Painlev\'e equation? 
\end{problem}
Since the solutions of the Painlev\'e equations are transcendental functions, it is difficult to make one in the same way as linear differential equations.

Recently, the author noticed that Professor P. Painlev\'e gave {\it Painlev\'e scheme} (see \cite{1,2} and \cite{Ince} P 323);
\begin{center}
\begin{tabular}{|c||c|} \hline 
Painlev\'e classification & Eq.\eqref{index1} \\ \hline
$k_i=1-\frac{1}{n_i}$ & $n_i$ \\ \hline
Type I:$\left(\frac{m+1}{m},\frac{m-1}{m} \right)$  & $(-m,m,1,1)$ \\ \hline
Type III:$\left(\frac{1}{2},\frac{1}{2},\frac{1}{2},\frac{1}{2} \right)$ & $(2,2,2,2)$ \\ \hline
Type IV:$\left(\frac{2}{3},\frac{2}{3},\frac{2}{3} \right)$ & $(1,3,3,3)$ \\ \hline
Type V:$\left(\frac{3}{4},\frac{3}{4},\frac{1}{2} \right)$ & $(1,4,4,2)$ \\ \hline
Type VI:$\left(\frac{5}{6},\frac{2}{3},\frac{1}{2} \right)$ & $(1,6,3,2)$ \\ \hline
\end{tabular}
\end{center}
In the Painlev\'e classification, Type II is omitted because it may be regarded as a degenerate case of Type III (see \cite{Ince} P 323). Here, we review the Painlev\'e exponent $k$ (see \cite{Ince} P 322);

\begin{align*}
\frac{d}{dz}\begin{pmatrix}
             W \\
             P 
             \end{pmatrix}&=\begin{pmatrix}
             P \\
             \frac{k P^2}{W} 
             \end{pmatrix}=\frac{1}{W} \begin{pmatrix}
             1 & 0  \\
             0 & k
             \end{pmatrix}\begin{pmatrix}
             W P \\
             P^2 
             \end{pmatrix}.
             \end{align*}
If $k \not= 1$, $W(z)=(A z+B)^{\frac{1}{1-k}} \ (A,B \in {\mathbb C})$. In the case of $k=1$, $W(z)=e^{A z+B}$.

Setting $\frac{1}{1-k}=n \ (n \in {\mathbb Z})$, we see the following relation:
\begin{center}
\begin{tabular}{|c||c|} \hline 
Painlev\'e exponent & (continued) ratio of local index (Definition; see \eqref{b}) \\ \hline
$k=1-\frac{1}{n}$ & $n$ \\ \hline
\end{tabular}
\end{center}
We remark that $n \in {\mathbb Z}$ for necessary condition of Painlev\'e property. We also note that in \cite{Ince} Page 322 $\frac{1}{1+k}=n \ ( k=1+\frac{1}{n})$ were given. 

In this note, we review the notion of Painlev\'e scheme of the sixth Painlev\'e equation from the viewpoint of accessible singular point and its local index in the Hirzebruch surface of degree two ${\Sigma_2}$.

At first, we consider the case of the sixth Painlev\'e equation. The sixth Painlev\'e equation is equivalent to the following Hamiltonian system (see \cite{T1,O3}):
\begin{equation}\label{PVI}
  \left\{
  \begin{aligned}
   \frac{dx}{dt} &=\frac{\partial H_{VI}}{\partial y}=\frac{1}{t(t-1)}\{2y(x-t)(x-1)x-(\alpha_0-1)(x-1)x-\alpha_3(x-t)x\\
                 & \qquad -\alpha_4(x-t)(x-1)\},\\
   \frac{dy}{dt} &=-\frac{\partial H_{VI}}{\partial x}=\frac{1}{t(t-1)}[-\{(x-t)(x-1)+(x-t)x+(x-1)x\}y^2+\{(\alpha_0-1)(2x-1)\\
                 & \qquad +\alpha_3(2x-t)+\alpha_4(2x-t-1)\}y-\alpha_2(\alpha_1+\alpha_2)]
   \end{aligned}
  \right. 
\end{equation}
with the polynomial Hamiltonian
\begin{align}\label{HVI}
\begin{split}
&H_{VI}(x,y,t;\alpha_0,\alpha_1,\alpha_2,\alpha_3,\alpha_4)\\
&=\frac{1}{t(t-1)}[y^2(x-t)(x-1)x-\{(\alpha_0-1)(x-1)x+\alpha_3(x-t)x\\
&+\alpha_4(x-t)(x-1)\}y+\alpha_2(\alpha_1+\alpha_2)x]  \quad (\alpha_0+\alpha_1+2\alpha_2+\alpha_3+\alpha_4=1). 
\end{split}
\end{align}
Since each right hand side of this system is polynomial with respect to $x,y$, by Cauchy's existence and uniqueness theorem of solutions, there exists unique holomorphic solution with initial values $(x,y)=(x_0,y_0) \in {\mathbb C}^2$.

Let us extend the regular vector field defined on ${\mathbb C}^2 \times B$
$$
v=\frac{\partial}{\partial t}+\frac{\partial H_{VI}}{\partial y}\frac{\partial}{\partial x}-\frac{\partial H_{VI}}{\partial x}\frac{\partial}{\partial y}
$$
to a rational vector field on ${\Sigma_2} \times B$, where $B={\mathbb C}-\{0,1\}$.

\begin{figure}
\unitlength 0.1in
\begin{picture}( 24.8000, 16.3000)( 19.9000,-20.0000)
%
\special{pn 8}%
\special{pa 1990 590}%
\special{pa 4320 590}%
\special{fp}%
%
\special{pn 8}%
\special{pa 2010 1810}%
\special{pa 4310 1810}%
\special{fp}%
%
\special{pn 8}%
\special{pa 2410 410}%
\special{pa 2410 2000}%
\special{fp}%
%
\special{pn 8}%
\special{pa 3990 410}%
\special{pa 3990 1990}%
\special{fp}%
\put(44.7000,-13.0000){\makebox(0,0)[lb]{${\Sigma}_{2}$}}%
%
\special{pn 20}%
\special{pa 2410 590}%
\special{pa 2710 590}%
\special{fp}%
\special{sh 1}%
\special{pa 2710 590}%
\special{pa 2644 570}%
\special{pa 2658 590}%
\special{pa 2644 610}%
\special{pa 2710 590}%
\special{fp}%
%
\special{pn 20}%
\special{pa 2410 590}%
\special{pa 2410 870}%
\special{fp}%
\special{sh 1}%
\special{pa 2410 870}%
\special{pa 2430 804}%
\special{pa 2410 818}%
\special{pa 2390 804}%
\special{pa 2410 870}%
\special{fp}%
\put(24.5000,-5.4000){\makebox(0,0)[lb]{$q$}}%
\put(21.4000,-8.3000){\makebox(0,0)[lb]{$p$}}%
\put(44.7000,-18.7000){\makebox(0,0)[lb]{$D^{(0)} \cong {\mathbb P}^1$}}%
%
\special{pn 20}%
\special{pa 2410 1810}%
\special{pa 2410 1540}%
\special{fp}%
\special{sh 1}%
\special{pa 2410 1540}%
\special{pa 2390 1608}%
\special{pa 2410 1594}%
\special{pa 2430 1608}%
\special{pa 2410 1540}%
\special{fp}%
%
\special{pn 20}%
\special{pa 2400 1810}%
\special{pa 2680 1810}%
\special{fp}%
\special{sh 1}%
\special{pa 2680 1810}%
\special{pa 2614 1790}%
\special{pa 2628 1810}%
\special{pa 2614 1830}%
\special{pa 2680 1810}%
\special{fp}%
%
\special{pn 20}%
\special{pa 3990 590}%
\special{pa 3990 870}%
\special{fp}%
\special{sh 1}%
\special{pa 3990 870}%
\special{pa 4010 804}%
\special{pa 3990 818}%
\special{pa 3970 804}%
\special{pa 3990 870}%
\special{fp}%
%
\special{pn 20}%
\special{pa 3980 590}%
\special{pa 3710 590}%
\special{fp}%
\special{sh 1}%
\special{pa 3710 590}%
\special{pa 3778 610}%
\special{pa 3764 590}%
\special{pa 3778 570}%
\special{pa 3710 590}%
\special{fp}%
%
\special{pn 20}%
\special{pa 3980 1810}%
\special{pa 3980 1530}%
\special{fp}%
\special{sh 1}%
\special{pa 3980 1530}%
\special{pa 3960 1598}%
\special{pa 3980 1584}%
\special{pa 4000 1598}%
\special{pa 3980 1530}%
\special{fp}%
%
\special{pn 20}%
\special{pa 3970 1800}%
\special{pa 3720 1800}%
\special{fp}%
\special{sh 1}%
\special{pa 3720 1800}%
\special{pa 3788 1820}%
\special{pa 3774 1800}%
\special{pa 3788 1780}%
\special{pa 3720 1800}%
\special{fp}%
\put(36.3000,-5.4000){\makebox(0,0)[lb]{$z_1$}}%
\put(40.8000,-8.2000){\makebox(0,0)[lb]{$w_1$}}%
\put(20.9000,-17.3000){\makebox(0,0)[lb]{$w_2$}}%
\put(24.6000,-20.0000){\makebox(0,0)[lb]{$z_2$}}%
\put(40.8000,-17.3000){\makebox(0,0)[lb]{$w_3$}}%
\put(36.7000,-20.0000){\makebox(0,0)[lb]{$z_3$}}%
\end{picture}%
\label{fig:E6figure1}
\caption{Hirzebruch surface ${\Sigma_2}$}
\end{figure}
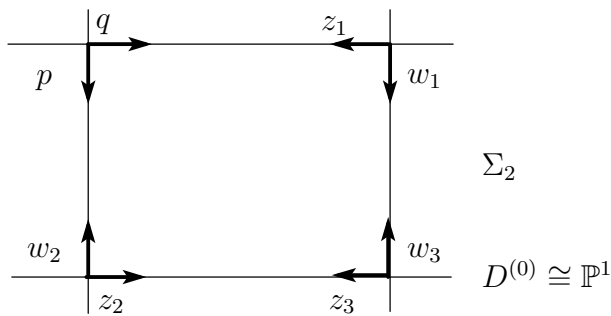

Here, we review the Hirzebruch surface ${\Sigma_2}$, which is obtained by gluing four copies of ${\mathbb C}^2$ via the following identification:
\begin{align}
\begin{split}
&U_j \cong {\mathbb C}^2 \ni (z_j,w_j) \ (j=0,1,2,3)\\
&z_0=x, \ w_0=y, \quad z_1=\frac{1}{x}, \ w_1=-(xy+\alpha_2)x,\\
&z_2=z_0, \ w_2=\frac{1}{w_0}, \quad z_3=z_1, \ w_3=\frac{1}{w_1}.
\end{split}
\end{align}
We define a divisor $D^{(0)}$ on ${\Sigma_2}$:
\begin{equation}
D^{(0)}=\{(z_2,w_2) \in U_2|w_2=0\} \cup \{(z_3,w_3) \in U_3|w_3=0\} \cong {\mathbb P}^1.
\end{equation}
The self-intersection number of $D^{(0)}$ is given by
\begin{equation}
(D^{(0)})^2=2.
\end{equation}
In the coordinate system $(z_1,w_1)$ the right hand side of this system is polynomial with respect to $z_1,w_1$. This compactification was found by Professor K. Okamoto (see \cite{O3}).

However, on the boundary divisor $D^{(0)} \cong {\mathbb P}^1$ this system has a pole in each coordinate system $(z_i,w_i) \ i=2,3$, whose order is one.
 By calculating its accessible singular points on $D^{(0)}$, we can obtain simple four singular points $z_2=0,1,t,\infty$ (see Definition \ref{Def1}).

By resolving all singular points, we can construct the space of initial conditions of the Painlev\'e VI system (see \cite{O3}). This space parametrizes all meromorphic solutions including holomorphic solutions.

 Conversely, we can recover the Painlev\'e VI system by all patching data of its space of initial conditions. In this note, we decompose its patching data into the pair of accessible singular point and local index $(n_i,1), \ n_i \in {\mathbb C}$ around each singular point $z_2=c_i \ (c_i \in \{0,1,t,\infty \})$ in the Hirzebruch surface of degree two ${\Sigma_2}$;
\begin{equation}
\begin{pmatrix}
z_2=0 & z_2=1 & z_2=t & z_2=\infty\\
\begin{pmatrix}
n_1 & \alpha_4\\
0 & 1 
\end{pmatrix} & \begin{pmatrix}
n_2 & \alpha_3\\
0 & 1 
\end{pmatrix} & \begin{pmatrix}
n_3 & \alpha_0\\
0 & 1 
\end{pmatrix} & \begin{pmatrix}
n_4 & \alpha_1\\
0 & 1 
\end{pmatrix}
\end{pmatrix},
\end{equation}
where the eigenvalues $n_i$ satisfy the following relation:
\begin{equation}
\frac{1}{n_1}+\frac{1}{n_2}+\frac{1}{n_3}+\frac{1}{n_4}=2, \quad (2 n_1 n_2 n_3 n_4-n_1 n_2 n_3-n_1 n_2 n_4-n_1 n_3 n_4-n_2 n_3 n_4=0).
\end{equation}
This equation has symmetry of symmetric group of degree four.

The key method is Painlev\'e $\alpha$-method for each accessible singular point. For example, let us consider the following differential system, which is equivalent to the Painlev\'e VI system \eqref{PVI} with \eqref{HVI};
{\Small
\begin{align}\label{11ZZZZ}
\begin{split}
&\frac{d}{dt}\begin{pmatrix}
             X \\
             Y 
             \end{pmatrix}=\frac{1}{t(t-1)Y} \{ t\begin{pmatrix}
             2 & -\alpha_4   \\
             0 & 1
             \end{pmatrix}\begin{pmatrix}
             X \\
             Y 
             \end{pmatrix}+\begin{pmatrix}
             -2(t+1) & \alpha_0-1+\alpha_4+t(\alpha_3+\alpha_4) \\
             0 & -2(t+1) 
             \end{pmatrix}\begin{pmatrix}
             X^2 \\
             XY 
             \end{pmatrix}\\
&+\begin{pmatrix}
             2 & -(\alpha_0+\alpha_3+\alpha_4)+1 \\
             0 & 3 
             \end{pmatrix}\begin{pmatrix}
             X^3 \\
             X^2 Y 
             \end{pmatrix} \}+\begin{pmatrix}
             0 \\
             \frac{-\{(\alpha_0-1)(2X-1)+\alpha_3(2X-t)+\alpha_4(2X-t-1)\}Y+\alpha_2(\alpha_1+\alpha_2)Y^2}{t(t-1)}
             \end{pmatrix},
             \end{split}
             \end{align}}where $(X,Y)=(x,1/y)$. This expansion is called {\it Painlev\'e expansion} (see \cite{Ince} P 322). We see that $(X,Y)=(0,0)$ is its accessible singular point. We remark that this system has a 1-parameter family of formal power series:
\begin{align}\label{Laurent series1}
\begin{split}
&X=-\frac{\alpha_4}{1-t_0}T+h T^2+{\mathcal O}(T^3), \quad Y=-\frac{1}{1-t_0}T+{\mathcal O}(T^2),
\end{split}
\end{align}
where $T:=t-t_0$, $h$ is its free parameter and the symbol ${\mathcal O}$ denotes Landau symbol. This formal power series coincides with known formal meromorphic solution (see \cite{Oka2}; P 212);
\begin{align}\label{Laurent series2}
\begin{split}
&x=-\frac{\alpha_4}{1-t_0}T+h T^2+{\mathcal O}(T^3), \quad y=-\frac{1-t_0}{T}[1+{\mathcal O}(T)], \quad (X,Y)=(x,1/y).
\end{split}
\end{align}

Now, let us make a change of variables $X,Y,t$ with a samll parameter $\alpha$:
\begin{equation}
X=\alpha Z, \quad Y=\alpha W, \quad t=t_0+\alpha T \quad (t_0 \in {\mathbb C}-\{0,1\}).
\end{equation}
Then the system can also be written in the new variables $Z,W,T$. This new system tends to the system as $\alpha \rightarrow 0$
\begin{align}\label{1ZZZZ}
\frac{d}{dT}\begin{pmatrix}
             Z \\
             W 
             \end{pmatrix}&=\frac{1}{W}\left\{\begin{pmatrix}
             \frac{2}{t_0-1} & -\frac{\alpha_4}{t_0-1}  \\
             0 & \frac{1}{t_0-1} 
             \end{pmatrix}\begin{pmatrix}
             Z \\
             W 
             \end{pmatrix}\right\}.
             \end{align}
We see that (continued) ratio of eigenvalues for the above matrix is given by $\frac{\frac{2}{t_0-1}}{\frac{1}{t_0-1} }=2$, which coincides with resonance data of formal power series \eqref{Laurent series1}.

Fixing $t=t_0$, this system is the system of the first order ordinary differential equation with constant coefficient. Let us solve this system explicitly;
\begin{equation}
Z(T)=C_2\{T+(t_0-1)C_1 \}^2+\frac{\alpha_4(T+(t_0-1)C_1)}{t_0-1}, \quad W(T)=\frac{T}{t_0-1}+C_1 \quad (C_1,C_2 \in {\mathbb C}).
\end{equation}
Thus, we can obtain single-valued solutions. For the Painlev\'e property, this is the necessary condition.

{\bf Painlev\'e VI case}
\begin{center}
\begin{tabular}{|c||c|} \hline
Equation & Painlev\'e VI system \eqref{PVI} with canonical Hamiltonian \eqref{HVI}  \\ \hline
Compactification & ${\Sigma_2}$: Hirzebruch surface of degree two  \\ \hline
Accessible singular points & $(z_2,w_2)=(0,0),(1,0),(t,0),(\infty,0)$  \\ \hline
Painlev\'e scheme & $
\begin{pmatrix}
z_2=0 & z_2=1 & z_2=t & z_2=\infty\\
\begin{pmatrix}
n_1 & \alpha_4\\
0 & 1 
\end{pmatrix} & \begin{pmatrix}
n_2 & \alpha_3\\
0 & 1 
\end{pmatrix} & \begin{pmatrix}
n_3 & \alpha_0\\
0 & 1 
\end{pmatrix} & \begin{pmatrix}
n_4 & \alpha_1\\
0 & 1 
\end{pmatrix}
\end{pmatrix}
$  \\ \hline
Relation of eigenvalues $n_i$ & $\frac{1}{n_1}+\frac{1}{n_2}+\frac{1}{n_3}+\frac{1}{n_4}=2$  \\ \hline
Painlev\'e VI case & $(n_1,n_2,n_3,n_4)=(2,2,2,2)$  \\ \hline
\end{tabular}
\end{center}
\begin{equation}
  \left\{
  \begin{aligned}
   \delta t(t-1)\frac{dx}{dt} =&n_1 n_2 n_3 x(x-1)(t-x)y+\{(2n_1n_2n_3-n_1n_2-n_1n_3-n_2n_3)\alpha_1-n_1n_2n_3\alpha_2\}x^2\\
&+\{-(2n_1n_2n_3-n_1n_2-n_1n_3-n_2n_3)\alpha_1+n_1n_2n_3\alpha_2+n_1 n_3 \alpha_3(t-1)\\
&+n_2n_3 \alpha_4 t \}x-n_2 n_3 \alpha_4 t,\\
\delta t(t-1)\frac{dy}{dt} =&[(n_1 n_2+n_2 n_3+n_1 n_3)x^2-\{n_1 n_2+n_2 n_3+(n_1+n_2)n_3 t\}x+n_2 n_3 t]y^2\\
&+[-\{2(2n_1n_2n_3-n_1n_2-n_1n_3-n_2n_3)\alpha_1\\
&+(-2n_1 n_2-2n_1 n_3-2n_2n_3+n_1n_2n_3)\alpha_2\}x\\
&+(-n_1 n_2-n_1 n_3-n_2n_3+2n_1n_2n_3)\alpha_1\\
&+\{(n_1 n_2-n_1-n_2)n_3 t-n_1 n_2-n_2 n_3\}\alpha_2\\
&-n_1 n_3 \alpha_3(t-1)-n_2 n_3\alpha_4 t]y-\alpha_2[(-n_1 n_2-n_1 n_3-n_2n_3+2n_1n_2n_3)\alpha_1\\
&+(-n_1 n_2-n_1 n_3-n_2n_3+n_1n_2n_3)\alpha_2].
   \end{aligned}
  \right. 
\end{equation}
Here, $\delta:=n_1 n_2 \alpha_0+(2n_1 n_2 n_3-n_1 n_2-n_1 n_3-n_2 n_3)\alpha_1-n_1 n_2 n_3 \alpha_2+n_1 n_3 \alpha_3+n_2 n_3 \alpha_4$.

This system is invariant under the following transformations\rm{:\rm} with the notation $(*)=(x,y,t;n_1,n_2,n_3,n_4;\alpha_0,\alpha_1,\ldots,\alpha_4)$  (See Section 4),
\begin{align*}
        s: (*) \rightarrow &(x+\frac{\alpha_2}{y},y,t;n_1,n_2,n_3,n_4;\alpha_0+\alpha_2-n_3 \alpha_2,\\
&\frac{(-n_1 n_2-n_1 n_3-n_2n_3+2n_1n_2n_3)\alpha_1+(-n_1 n_2-n_1 n_3-n_2n_3+n_1n_2n_3)\alpha_2}{-n_1 n_2-n_1 n_3-n_2n_3+2n_1n_2n_3},\\
&-\alpha_2,\alpha_3+\alpha_2-n_2 \alpha_2,\alpha_4+\alpha_2-n_1 \alpha_2), \\
        \pi_1: (*) \rightarrow &(1-x,-y,1-t;n_2,n_1,n_3,n_4;\alpha_0,\alpha_1,\alpha_2,\alpha_4,\alpha_3), \\
        \pi_2: (*) \rightarrow &\left(\frac{t-x}{t-1},-(t-1)y,\frac{t}{t-1};n_3,n_2,n_1,n_4;\alpha_4,\alpha_1,\alpha_2,\alpha_3,\alpha_0 \right), \\
        \pi_3: (*) \rightarrow &(\frac{1}{x},-(yx+\alpha_2)x,\frac{1}{t};n_4,n_2,n_3,n_1;\alpha_0,\frac{\alpha_4 n_2 n_3 n_4}{n_1(2n_1n_2n_3-n_1 n_2-n_1 n_3-n_2n_3)},\\
&\alpha_2,\alpha_3,\frac{\alpha_1 (2n_1n_2n_3-n_1 n_2-n_1 n_3-n_2n_3)}{n_1 n_2 n_3}).
        \end{align*}

\section{Review of accessible singularity and local index}
Let us review the notion of {\it accessible singularity}. Let $B$ be a connected open domain in $\mathbb C$ and $\pi : {\mathcal W} \longrightarrow B$ a smooth proper holomorphic map. We assume that ${\mathcal H} \subset {\mathcal W}$ is a normal crossing divisor which is flat over $B$. Let us consider a rational vector field $\tilde v$ on $\mathcal W$ satisfying the condition
\begin{equation*}
\tilde v \in H^0({\mathcal W},\Theta_{\mathcal W}(-\log{\mathcal H})({\mathcal H})).
\end{equation*}
Fixing $t_0 \in B$ and $P \in {\mathcal W}_{t_0}$, we can take a local coordinate system $(x_1,\ldots ,x_n)$ of ${\mathcal W}_{t_0}$ centered at $P$ such that ${\mathcal H}_{\rm smooth \rm}$ can be defined by the local equation $x_1=0$.
Since $\tilde v \in H^0({\mathcal W},\Theta_{\mathcal W}(-\log{\mathcal H})({\mathcal H}))$, we can write down the vector field $\tilde v$ near $P=(0,\ldots ,0,t_0)$ as follows:
\begin{equation*}
\tilde v= \frac{\partial}{\partial t}+g_1 
\frac{\partial}{\partial x_1}+\frac{g_2}{x_1} 
\frac{\partial}{\partial x_2}+\cdots+\frac{g_n}{x_1} 
\frac{\partial}{\partial x_n}.
\end{equation*}
This vector field defines the following system of differential equations
\begin{equation}\label{39}
\frac{dx_1}{dt}=g_1(x_1,\ldots,x_n,t),\ \frac{dx_2}{dt}=\frac{g_2(x_1,\ldots,x_n,t)}{x_1},\cdots, \frac{dx_n}{dt}=\frac{g_n(x_1,\ldots,x_n,t)}{x_1}.
\end{equation}
Here $g_i(x_1,\ldots,x_n,t), \ i=1,2,\ldots ,n,$ are holomorphic functions defined near $P$.

\begin{definition}\label{Def1}
With the above notation, assume that the rational vector field $\tilde v$ on $\mathcal W$ satisfies the condition
$$
(A) \quad \tilde v \in H^0({\mathcal W},\Theta_{\mathcal W}(-\log{\mathcal H})({\mathcal H})).
$$
We say that $\tilde v$ has an {\it accessible singularity} at $P=(0,\dots ,0,t_0)$ if
\begin{equation}
\boxed{%
x_1=0 \ {\rm and \rm} \ g_i(0,\ldots,0,t_0)=0 \ {\rm for \rm} \ {\rm every \rm} \ i, \ 2 \leq i \leq n.
}%
\end{equation}
\end{definition}

If $P \in {\mathcal H}_{{\rm smooth \rm}}$ is not an accessible singularity, all solutions of the ordinary differential equation passing through $P$ are vertical solutions, that is, the solutions are contained in the fiber ${\mathcal W}_{t_0}$ over $t=t_0$. If $P \in {\mathcal H}_{\rm smooth \rm}$ is an accessible singularity, there may be a solution of \eqref{39} which passes through $P$ and goes into the interior ${\mathcal W}-{\mathcal H}$ of ${\mathcal W}$.

Here we review the notion of {\it local index}. Let $v$ be an algebraic vector field with an accessible singular point $\overrightarrow{p}=(0,\ldots,0)$ and $(x_1,\ldots,x_n)$ be a coordinate system in a neighborhood centered at $\overrightarrow{p}$. Assume that the system associated with $v$ near $\overrightarrow{p}$ can be written as

{\Small
\begin{align}\label{b}
\begin{split}
\frac{d}{dt}\begin{pmatrix}
             x_1 \\
             x_2 \\
             \vdots\\
             x_{n-1} \\
             x_n
             \end{pmatrix}=\frac{1}{x_1}\left\{\begin{bmatrix}
             a_{11} & 0 & 0 & \hdots & 0 \\
             a_{21} & a_{22} & 0 &  \hdots & 0 \\
             \vdots & \vdots & \ddots & 0 & 0 \\
             a_{(n-1)1} & a_{(n-1)2} & \hdots & a_{(n-1)(n-1)} & 0 \\
             a_{n1} & a_{n2} & \hdots & a_{n(n-1)} & a_{nn}
             \end{bmatrix}\begin{pmatrix}
             x_1 \\
             x_2 \\
             \vdots\\
             x_{n-1} \\
             x_n
             \end{pmatrix}+\begin{pmatrix}
             x_1h_1(x_1,\ldots,x_n,t) \\
             h_2(x_1,\ldots,x_n,t) \\
             \vdots\\
             h_{n-1}(x_1,\ldots,x_n,t) \\
             h_n(x_1,\ldots,x_n,t)
             \end{pmatrix}\right\},\\
              (h_i \in {\mathbb C}(t)[x_1,\ldots,x_n], \ a_{ij} \in {\mathbb C}(t))
             \end{split}
             \end{align}}
where $h_1$ is a polynomial which vanishes at $\overrightarrow{p}$ and $h_i$, $i=2,3,\ldots,n$ are polynomials of order at least 2 in $x_1,x_2,\ldots,x_n$, We call ordered set of the eigenvalues $(a_{11},a_{22},\cdots,a_{nn})$ {\it local index} at $\overrightarrow{p}$.

We are interested in the case with local index
\begin{equation}\label{integer}
\left(1,\frac{a_{22}(t)}{a_{11}(t)},\ldots,\frac{a_{nn}(t)}{a_{11}(t)} \right) \in {\mathbb Z}^{n}.
\end{equation}

If each component of $\left(1,\frac{a_{22}(t)}{a_{11}(t)},\ldots,\frac{a_{nn}(t)}{a_{11}(t)} \right)$ has the same sign, we may resolve the accessible singularity by blowing-up finitely many times. However, when different signs appear, we may need to both blow up and blow down.
\begin{center}
\begin{tabular}{|c||c|c|} \hline 
& $\left(\frac{a_{22}(t)}{a_{11}(t)},\ldots,\frac{a_{nn}(t)}{a_{11}(t)} \right)$  & Resolution of accessible sing. \\ \hline
Positive sign & ${\mathbb N}^{n-1}$ & Blowing-up \\ \hline
Different signs &  ${\mathbb Z}^{n-1}$ & both Blow up and Blow down  \\ \hline
\end{tabular}
\end{center}

The $\alpha$-test,
\begin{equation}\label{poiuy}
t=t_0+\alpha T, \quad x_i=\alpha X_i, \quad \alpha \rightarrow 0,
\end{equation}
yields the following reduced system:
\begin{align}\label{ppppppp}
\begin{split}
\frac{d}{dT}\begin{pmatrix}
             X_1 \\
             X_2 \\
             \vdots\\
             X_{n-1} \\
             X_n
             \end{pmatrix}=\frac{1}{X_1}\begin{bmatrix}
             a_{11}(t_0) & 0 & 0 & \hdots & 0 \\
             a_{21}(t_0) & a_{22}(t_0) & 0 &  \hdots & 0 \\
             \vdots & \vdots & \ddots & 0 & 0 \\
             a_{(n-1)1}(t_0) & a_{(n-1)2}(t_0) & \hdots & a_{(n-1)(n-1)}(t_0) & 0 \\
             a_{n1}(t_0) & a_{n2}(t_0) & \hdots & a_{n(n-1)}(t_0) & a_{nn}(t_0)
             \end{bmatrix}\begin{pmatrix}
             X_1 \\
             X_2 \\
             \vdots\\
             X_{n-1} \\
             X_n
             \end{pmatrix},
             \end{split}
             \end{align}
where $a_{ij}(t_0) \in {\mathbb C}$. Fixing $t=t_0$, this system is the system of the first order ordinary differential equation with constant coefficient. Let us solve this system. At first, we solve the first equation:
\begin{equation}
X_1(T)=a_{11}(t_0)T+C_1 \quad (C_1 \in {\mathbb C}).
\end{equation}
Substituting this into the second equation in \eqref{ppppppp}, we can obtain the first order linear ordinary differential equation:
\begin{equation}
\frac{dX_2}{dT}=\frac{a_{22}(t_0) X_2}{a_{11}(t_0)T+C_1}+a_{21}(t_0).
\end{equation}
In the case of $a_{11}(t_0) \not= a_{22}(t_0)$ we can solve explicitly:
\begin{equation}
X_2(T)=C_2(a_{11}(t_0)T+C_1)^{\frac{a_{22}(t_0)}{a_{11}(t_0)}}+\frac{a_{21}(t_0)(a_{11}(t_0)T+C_1)}{a_{11}(t_0)-a_{22}(t_0)} \quad (C_2 \in {\mathbb C}).
\end{equation}
This solution is a single-valued solution if and only if
$$
\frac{a_{22}(t_0)}{a_{11}(t_0)} \in {\mathbb Z}.
$$
In the case of $a_{11}(t_0)=a_{22}(t_0)$ we can solve explicitly;
\begin{equation}
X_2(T)=C_2(a_{11}(t_0)T+C_1)+\frac{a_{21}(t_0)(a_{11}(t_0)T+C_1){\rm Log}(a_{11}(t_0)T+C_1)}{a_{11}(t_0)} \quad (C_2 \in {\mathbb C}).
\end{equation}
This solution is a single-valued solution if and only if
\begin{equation}\label{eq:redlocal}
a_{21}(t_0)=0.
\end{equation}
Of course, $\frac{a_{22}(t_0)}{a_{11}(t_0)}=1 \in {\mathbb Z}$.
In the same way, we can obtain the solutions for each variables $(X_3,\ldots,X_n)$.

\begin{tabular}{|c|} \hline 
The conditions $\frac{a_{jj}(t)}{a_{11}(t)} \in {\mathbb Z}, \ (j=2,3,\ldots,n)$ are necessary condition in order to have\\
the Painlev\'e property. \\ \hline 
\end{tabular}

\begin{center}
\begin{tabular}{|c||c|c|} \hline 
& $\left(\frac{a_{22}(t)}{a_{11}(t)},\ldots,\frac{a_{nn}(t)}{a_{11}(t)} \right)$  & Movable singularities \\ \hline
Painlev\'e type & ${\mathbb Z}$ & Only pole  \\ \hline
Other Non-Linear Equation &  ${\mathbb Q},{\mathbb R}$ and ${\mathbb C}$ & Algebraic sing. or others  \\ \hline
\end{tabular}
\end{center}

For example, we consider the Painlev\'e VI equation. Let us calculate its accessible singularities.

In the coordinate system $(X,Y)=(z_2,w_2)=(x.1/y)$ we can rewrite the system given by
\begin{equation}
  \left\{
  \begin{aligned}
   \frac{dX}{dt} &= \frac{2 X (X - 1) (X - t) }{t(t - 
     1)  Y} + \frac{(1 - \alpha_0 - \alpha_3 - \alpha_4) 
    X^2 + (\alpha_0 + \alpha_4 - 1 + (\alpha_3 + \alpha_4) t ) 
    X - \alpha_4 t }{t(t - 1) },\\
   \frac{dY}{dt} &=\frac{1}{t - 1}  + \frac{(3 X - 2 (t + 1)) X + 
  \alpha_2 (\alpha_1 + \alpha_2) Y^2 }{t(t - 1)}\\
& + \frac{ 
  (\alpha_0 + \alpha_4 - 1 + (\alpha_3 + \alpha_4) t - 
     2 (\alpha_0 + \alpha_3 + \alpha_4 - 1) X)Y}{t(t - 1)}
.
   \end{aligned}
  \right. 
\end{equation}
By a direct calculation, we can obtain some {\bf accessible singular points};

$\{(X,Y)|Y=0, \ X (X - 1) (X - t) =0 \}=\{(X,Y)= (0,0),(1,0),(t,0) \}$.

Next, let us rewrite the system centered at each singular point $X=0,1,t,\infty$.

1. By taking the coordinate system $(X,Y)=(z_2,w_2)$ centered at the point $(z_2,w_2)=(0,0)$, the system is given by

{\Small
\begin{align*}
&\frac{d}{dt}\begin{pmatrix}
             X \\
             Y 
             \end{pmatrix}=\frac{1}{t(t-1)Y} \{ t\begin{pmatrix}
             2 & -\alpha_4   \\
             0 & 1
             \end{pmatrix}\begin{pmatrix}
             X \\
             Y 
             \end{pmatrix}+\begin{pmatrix}
             -2(t+1) & \alpha_0-1+\alpha_4+t(\alpha_3+\alpha_4) \\
             0 & -2(t+1) 
             \end{pmatrix}\begin{pmatrix}
             X^2 \\
             XY 
             \end{pmatrix}\\
&+\begin{pmatrix}
             2 & -(\alpha_0+\alpha_3+\alpha_4)+1 \\
             0 & 3 
             \end{pmatrix}\begin{pmatrix}
             X^3 \\
             X^2 Y 
             \end{pmatrix} \}+\begin{pmatrix}
             0 \\
             \frac{-\{(\alpha_0-1)(2X-1)+\alpha_3(2X-t)+\alpha_4(2X-t-1)\}Y+\alpha_2(\alpha_1+\alpha_2)Y^2}{t(t-1)}
             \end{pmatrix}.
             \end{align*}}
This expansion is called {\it Painlev\'e expansion}.  

Now, let us make a change of variables $X,Y,t$ with a samll parameter $\alpha$:
\begin{equation}
X=\alpha Z, \quad Y=\alpha W, \quad t=t_0+\alpha T \quad (t_0 \in {\mathbb C}-\{0,1\}).
\end{equation}
Then the system can also be written in the new variables $Z,W,T$. This new system tends to the system as $\alpha \rightarrow 0$
\begin{align}\label{ZZZZ}
\frac{d}{dT}\begin{pmatrix}
             Z \\
             W 
             \end{pmatrix}&=\frac{1}{W}\left\{\begin{pmatrix}
             \frac{2}{t_0-1} & -\frac{\alpha_4}{t_0-1}  \\
             0 & \frac{1}{t_0-1} 
             \end{pmatrix}\begin{pmatrix}
             Z \\
             W 
             \end{pmatrix}\right\}.
             \end{align}
Fixing $t=t_0$, this system is the system of the first order ordinary differential equation with constant coefficient. Let us solve this system. At first, we solve the second equation:
\begin{equation}
W(T)=\frac{T}{t_0-1}+C_1 \quad (C_1 \in {\mathbb C}).
\end{equation}
Substituting this into the first equation in \eqref{ZZZZ}, we can obtain the first order linear ordinary differential equation:
\begin{equation}
\frac{dZ}{dT}=\frac{t_0-1}{T+C_1(t_0-1)}\left(\frac{2}{t_0-1}Z-\frac{\alpha_4}{t_0-1} \left(\frac{T}{t_0-1}+C_1 \right) \right).
\end{equation}
We can solve explicitly:
\begin{equation}
Z(T)=C_2\{T+(t_0-1)C_1 \}^2+\frac{\alpha_4(T+(t_0-1)C_1)}{t_0-1} \quad (C_2 \in {\mathbb C}).
\end{equation}
Thus, we can obtain single-valued solutions. For the Painlev\'e property, this is the necessary condition.

\vspace{0.5cm}
In the same way, we can obtain the following:

2. By taking the coordinate system $(X,Y)=(z_2-1,w_2)$ centered at the point $(z_2,w_2)=(1,0)$, the system is given by
\begin{align*}
\frac{d}{dt}\begin{pmatrix}
             X \\
             Y 
             \end{pmatrix}&=\frac{1}{Y}\left\{\begin{pmatrix}
             -\frac{2}{t} & \frac{\alpha_3}{t}  \\
             0 & -\frac{1}{t} 
             \end{pmatrix}\begin{pmatrix}
             X \\
             Y 
             \end{pmatrix}+\cdots\right\}.
             \end{align*}

3. By taking the coordinate system $(X,Y)=(z_2-t,w_2)$ centered at the point $(z_2,w_2)=(t,0)$, the system is given by
\begin{align*}
\frac{d}{dt}\begin{pmatrix}
             X \\
             Y 
             \end{pmatrix}&=\frac{1}{Y}\left\{\begin{pmatrix}
             2 & -\alpha_0  \\
             0 & 1 
             \end{pmatrix}\begin{pmatrix}
             X \\
             Y 
             \end{pmatrix}+\cdots\right\}.
             \end{align*}

4. By taking the coordinate system $(X,Y)=(z_3,w_3)$ centered at the point $(z_3,w_3)=(0,0)$, the system is given by
\begin{align*}
\frac{d}{dt}\begin{pmatrix}
             X \\
             Y 
             \end{pmatrix}&=\frac{1}{Y}\left\{\begin{pmatrix}
             \frac{2}{t(t-1)} & -\frac{\alpha_1}{t(t-1)}  \\
             0 & \frac{1}{t(t-1)} 
             \end{pmatrix}\begin{pmatrix}
             X \\
             Y 
             \end{pmatrix}+\cdots\right\}.
             \end{align*}
Thus, we have proved that the Painlev\'e VI Hamiltonian system \eqref{PVI},\eqref{HVI} passes the Painlev\'e $\alpha$-test for all accessible singular points $X=0,1,t,\infty;$

\begin{equation*}
\begin{pmatrix}
X=0 & X=1 & X=t & X=\infty\\
\frac{1}{t-1}\begin{pmatrix}
2 & -\alpha_4\\
0 & 1 
\end{pmatrix} & -\frac{1}{t}\begin{pmatrix}
2 & -\alpha_3\\
0 & 1 
\end{pmatrix} & \begin{pmatrix}
2 & -\alpha_0\\
0 & 1 
\end{pmatrix} & \frac{1}{t(t-1)}\begin{pmatrix}
2 & -\alpha_1\\
0 & 1 
\end{pmatrix}
\end{pmatrix}.
\end{equation*}

\begin{center}
\begin{tabular}{|c|} \hline 
The pair of accessible singular points and matrix of linear approximation\\
around each point is called {\it Painlev\'e scheme} (see \cite{Ince} P323, cf. \cite{1,2}). \\ \hline 
\end{tabular}
\end{center}

\section{Recovery of the Painlev\'e VI system}

Let us consider the system of the first order ordinary differential equations of polynomial type.
\begin{equation}\label{1}
  \left\{
  \begin{aligned}
   \frac{dx}{dt} &=f_1(x,y),\\
   \frac{dy}{dt} &=f_2(x,y) \quad (f_i \in {\mathbb C}(t)[x,y]).
   \end{aligned}
  \right. 
\end{equation}
We assume that associated vector field defined on ${\mathbb C}^2 \times B$
$$
v=\frac{\partial}{\partial t}+f_1(x,y)\frac{\partial}{\partial x}+f_2(x,y)\frac{\partial}{\partial y}
$$
belongs in 
$$
v \in H^0({\Sigma_2},\Theta_{\Sigma_2}(-\log{\mathcal H})({\mathcal H})).
$$
This condition is equivalent to the following:
\begin{enumerate}\label{2}
\item Holomorphy in the coordinate system $(x_1,y_1)=(1/x,-(xy+\alpha_2)x)$,
\item In the coordinate system $(X,Y)=(x,1/y)$, the differential system \eqref{1} must be taken of the form:
\begin{equation}
  \left\{
  \begin{aligned}
   \frac{dX}{dt} &=\frac{F_1(X,Y)}{Y},\\
   \frac{dY}{dt} &=F_2(X,Y) \quad (F_i \in {\mathbb C}(t)[X,Y]).
   \end{aligned}
  \right. 
\end{equation}
\end{enumerate}

\begin{proposition}\label{Proposition3.1}
Under above assumptions 1 and 2, the system \eqref{1} is given by
\begin{equation}\label{4}
  \left\{
  \begin{aligned}
   \frac{dx}{dt} &=a_1x^3y+a_2x^2y+a_5xy+a_7y+\frac{1}{2}((3a_1+2a_3)\alpha_2-a_4)x^2+((a_2+a_9)\alpha_2-a_6)x+a_8,\\
   \frac{dy}{dt} &=a_3x^2y^2+a_9xy^2+a_{10}y^2+a_4xy+a_6y+\frac{1}{2}((a_1 \alpha_2+a_4)\alpha_2 \quad (a_i \in {\mathbb C}(t)).
   \end{aligned}
  \right. 
\end{equation}
Here, $a_i=a_i(t), \ (i=1,2,\ldots,10)$ are undetermined coefficients.
\end{proposition}

In the coordinate system $(x_1,y_1)=(1/x,-(xy+\alpha_2)x)$, the system \eqref{4} can be rewritten as follows:
\begin{equation}\label{444}
  \left\{
  \begin{aligned}
   \frac{dx_1}{dt} =&a_7 x_1^4y_1+a_5 x_1^3 y_1+a_2 x_1^2y_1+a_1x_1y_1+\alpha_2 a_7x_1^3+(\alpha_2 a_5-a_8)x_1^2+(a_6-\alpha_2 a_9)x_1\\
&-\frac{1}{2}\alpha_2(a_1+2a_3)+\frac{a_4}{2},\\
   \frac{dy_1}{dt} =&-2a_7x_1^3y_1^2-(2a_5+a_{10})x_1^2y_1^2-3\alpha_2a_7x_1^2y_1-(2a_2+a_9)x_1y_1^2-(2a_1+a_3)y_1^2\\
&-\{\alpha_2(3a_5+2a_{10})-2a_8 \}x_1y_1-\alpha_2^2 a_7x_1-(\alpha_2 a_2+a_6)y_1-\alpha_2 \{\alpha_2(a_5+a_{10})-a_8 \}.
   \end{aligned}
  \right. 
\end{equation}

{\bf Proof of Proposition \ref{Proposition3.1}.}

(i) Degree of polynomials $f_i(x,y)$ with respect to $y$

If the system \eqref{1} belongs in $H^0({\Sigma_2},\Theta_{\Sigma_2}(-\log{\mathcal H})({\mathcal H}))$,
in the coordinate system $(X,Y)=(x,1/y)$ this system must be taken of the form:
\begin{equation}\label{3}
  \left\{
  \begin{aligned}
   \frac{dX}{dt} &=\frac{F_1(X,Y)}{Y},\\
   \frac{dY}{dt} &=F_2(X,Y) \quad (F_i \in {\mathbb C}(t)[X,Y]).
   \end{aligned}
  \right. 
\end{equation}
By this condition, we see that the system \eqref{1} must be taken of the form:
\begin{equation}\label{111111}
  \left\{
  \begin{aligned}
   \frac{dx}{dt} &=b_1(x)+b_2(x)y,\\
   \frac{dy}{dt} &=b_3(x)+b_4(x)y+b_5(x)y^2 \quad (b_i \in {\mathbb C}(t)[x]).
   \end{aligned}
  \right. 
\end{equation}
Here, the degree of each $b_i$  with respect to $x$ is given by
\begin{equation}
deg(b_1)=l, \ deg(b_2)=m, \ deg(b_3)=n, \ deg(b_4)=p, \ deg(b_5)=r,
\end{equation}
where $l,m,n,p,r \in {\mathbb N}$.

(ii) Holomorphy in the coordinate system $(x_1,y_1)=(1/x,-(xy+\alpha_2)x)$

In the coordinate system $(x_1,y_1)=(1/x,-(xy+\alpha_2)x)$, the first equation of the system \eqref{111111} is given by
\begin{equation}
\frac{d x_1}{dt}=-x_1^2\left\{b_1\left(\frac{1}{x_1} \right)+b_2\left(\frac{1}{x_1} \right)(-x_1^2 y_1-\alpha_2 x_1) \right\}.
\end{equation}
Since the right hand side of this system must be polynomial with respect to $x_1$, we compare two terms
\begin{equation}
  \left\{
  \begin{aligned}
   b_1\left(\frac{1}{x_1} \right) &=\frac{b_1^{(l)}}{x_1^l}+\cdots,\\
   -\alpha_2 x_1 b_2\left(\frac{1}{x_1} \right) &=-\alpha_2 \frac{b_2^{(m)}}{x_1^{m-1}}+\cdots.
   \end{aligned}
  \right. 
\end{equation}
Since $b_1^{(l)} \not=0$ and $b_2^{(m)} \not=0$, we can obtain
\begin{equation}
l=m-1
\end{equation}
Next, we compare the term involving $y_1$:
\begin{equation}
-x_1^4 b_2\left(\frac{1}{x_1} \right)y_1=-x_1^4 \left(\frac{b_2^{(m)}}{x_1^m}+\frac{b_2^{(m-1)}}{x_1^{m-1}}+\cdots \right)y_1 \quad (b_2^{(j)} \in {\mathbb C}(t)).
\end{equation}
If this becomes polynomial with respect to $x_1,y_1$, 
\begin{equation}
m=4
\end{equation}
In the same way, we can obtain
\begin{equation}
deg(b_1)=3, \ deg(b_2)=4, \ deg(b_3)=1, \ deg(b_4)=2, \ deg(b_5)=3.
\end{equation}
Finally, by comparing undetermined coefficients, we can obtain the conclusion.

\vspace{0.5cm}
For the system \eqref{4}, by giving the following Painlev\'e scheme we can recover the Painlev\'e VI system with the polynomial Hamiltonian $H_{VI}$.

\begin{theorem}\label{Theorem3.1}
For the system \eqref{4}, we give the following Painlev\'e scheme:
\begin{equation}\label{asd}
\begin{pmatrix}
X=0 & X=1 & X=t & X=\infty\\
f_0\begin{pmatrix}
2 & -\alpha_4\\
0 & 1 
\end{pmatrix} & f_1\begin{pmatrix}
2 & -\alpha_3\\
0 & 1 
\end{pmatrix} & f_2\begin{pmatrix}
2 & -\alpha_0\\
0 & 1 
\end{pmatrix} & f_3\begin{pmatrix}
2 & -\alpha_1\\
0 & 1 
\end{pmatrix}
\end{pmatrix}.
\end{equation}
Here, $X=0,1,t,\infty$ are accessible singular points, $f_i \in {\mathbb C}(t)$ and $\alpha_i$ are constant parameters.
Then, this system coincides with the Painlev\'e VI system with the polynomial Hamiltonian $H_{VI}$.
\end{theorem}

{\bf Proof of Theorem \ref{Theorem3.1}.} \quad At first, we can rewrite the system \eqref{4} in the coordinate system $(X,Y)=(x,1/y)$ centered at $(X,Y)=(0,0)$
\begin{equation}\label{5}
  \left\{
  \begin{aligned}
   \frac{dX}{dt} &=\frac{a_1X^3+a_2X^2+a_5X+a_7}{Y}+\frac{1}{2}\{(3a_1+2a_3)\alpha_2-a_4\}X^2+\{(a_2+a_9)\alpha_2-a_6\}X+a_8,\\
   \frac{dY}{dt} &=-a_{10}-a_9X-a_3X^2-a_4XY-a_6Y-\frac{1}{2}(a_1 \alpha_2+a_4)\alpha_2Y^2 \quad (a_i \in {\mathbb C}(t)).
   \end{aligned}
  \right. 
\end{equation}
By Definition \ref{Def1}, we can calculate the accessible singular points
\begin{equation}\label{6}
Y=0, \quad a_1X^3+a_2X^2+a_5X+a_7=0.
\end{equation}
By the assumption, $X=Y=0$ is a solution of the system \eqref{6}. Thus, we obtain the condition
$$
a_7=0.
$$

By the assumption, the matrix of linear approximation around $X=0$ is given by
\begin{equation}
f_0\begin{pmatrix}
2 & -\alpha_4\\
0 & 1 
\end{pmatrix}.
\end{equation}
So, we obtain
$$
a_{5}=-2a_{10}, \quad a_{8}=\alpha_4 a_{10}.
$$

In the same way, we can obtain the conditions at each singular point $X=1,t,\infty$. Thus, we have completed the proof of Theorem \ref{Theorem3.1}.

\section{A generalization of the Painlev\'e VI system}
By generalizing the eigenvalues $(2,1)$ to $(n_i,1)$ $(i=1,2,3,4)$ in \eqref{asd}, we will construct a generalization of the Painlev\'e VI system.

\begin{theorem}\label{Theorem4.1}
For the system \eqref{4}, we give the following Painlev\'e scheme:
\begin{equation}\label{555}
\begin{pmatrix}
X=0 & X=1 & X=t & X=\infty\\
f_0\begin{pmatrix}
n_1 & \alpha_4\\
0 & 1 
\end{pmatrix} & f_1\begin{pmatrix}
n_2 & \alpha_3\\
0 & 1 
\end{pmatrix} & f_2\begin{pmatrix}
n_3 & \alpha_0\\
0 & 1 
\end{pmatrix} & f_3\begin{pmatrix}
n_4 & \alpha_1\\
0 & 1 
\end{pmatrix}
\end{pmatrix}.
\end{equation}
Here, $X=0,1,t,\infty$ are accessible singular points, $f_i \in {\mathbb C}(t)$, $n_i \in {\mathbb C}$ and $\alpha_i$ are constant parameters.
Then, this system coincides with
\begin{equation}\label{erty}
  \left\{
  \begin{aligned}
   \delta t(t-1)\frac{dx}{dt} =&n_1 n_2 n_3 x(x-1)(t-x)y+\{(2n_1n_2n_3-n_1n_2-n_1n_3-n_2n_3)\alpha_1-n_1n_2n_3\alpha_2\}x^2\\
&+\{-(2n_1n_2n_3-n_1n_2-n_1n_3-n_2n_3)\alpha_1+n_1n_2n_3\alpha_2+n_1 n_3 \alpha_3(t-1)\\
&+n_2n_3 \alpha_4 t \}x-n_2 n_3 \alpha_4 t,\\
\delta t(t-1)\frac{dy}{dt} =&[(n_1 n_2+n_2 n_3+n_1 n_3)x^2-\{n_1 n_2+n_2 n_3+(n_1+n_2)n_3 t\}x+n_2 n_3 t]y^2\\
&+[-\{2(2n_1n_2n_3-n_1n_2-n_1n_3-n_2n_3)\alpha_1\\
&+(-2n_1 n_2-2n_1 n_3-2n_2n_3+n_1n_2n_3)\alpha_2\}x\\
&+(-n_1 n_2-n_1 n_3-n_2n_3+2n_1n_2n_3)\alpha_1\\
&+\{(n_1 n_2-n_1-n_2)n_3 t-n_1 n_2-n_2 n_3\}\alpha_2\\
&-n_1 n_3 \alpha_3(t-1)-n_2 n_3\alpha_4 t]y-\alpha_2[(-n_1 n_2-n_1 n_3-n_2n_3+2n_1n_2n_3)\alpha_1\\
&+(-n_1 n_2-n_1 n_3-n_2n_3+n_1n_2n_3)\alpha_2].
   \end{aligned}
  \right. 
\end{equation}
Here, $\delta:=n_1 n_2 \alpha_0+(2n_1 n_2 n_3-n_1 n_2-n_1 n_3-n_2 n_3)\alpha_1-n_1 n_2 n_3 \alpha_2+n_1 n_3 \alpha_3+n_2 n_3 \alpha_4$.
\end{theorem}

\begin{proposition}
By using the Painlev\'e $\alpha$-method, we see that the system \eqref{erty} has movable branch points.
\end{proposition}

\begin{proposition}\label{Proposition4.1}
The eigenvalues $n_i$ satisfy the following relation:
\begin{equation}\label{234}
\frac{1}{n_1}+\frac{1}{n_2}+\frac{1}{n_3}+\frac{1}{n_4}=2.
\end{equation}
\end{proposition}

{\bf Proof of Proposition \ref{Proposition4.1}} \quad For the system \eqref{5}, we put $f:=a_1X^3+a_2X^2+a_5X+a_7$. Since the cubic equation $f=0$ has the solutions $X=0,1,t$,  from the relation between solution and coefficient we obtain
\begin{align*}
-\frac{a_2}{a_1}&=0+1+t=t+1,\\
\frac{a_5}{a_1}&=0 \times 1+0 \times t+1 \times t=t.
\end{align*}
We summarize that
\begin{align*}
a_2&=-(t+1)a_1,\\
a_5&=t a_1.
\end{align*}
The equation $f$ is given by
\begin{align*}
f&=a_1(X^2-(t+1)X+t)X\\
&=a_1(X-0)(X-1)(X-t).
\end{align*}
Thus, we can obtain
\begin{equation}
  \left\{
  \begin{aligned}
   \frac{dX}{dt} &=\frac{a_1(X-0)(X-1)(X-t)}{Y}+\frac{1}{2}\{(3a_1+2a_3)\alpha_2-a_4\}X^2+\{(a_2+a_9)\alpha_2-a_6\}X+a_8,\\
   \frac{dY}{dt} &=-a_{10}-a_9X-a_3X^2-a_4XY-a_6Y-\frac{1}{2}(a_1 \alpha_2+a_4)\alpha_2Y^2 \quad (a_i \in {\mathbb C}(t)).
   \end{aligned}
  \right. 
\end{equation}
Next, by giving the eigenvalues of the matrix of linear approximation around each point $X=0,1,t$, we can obtain
\begin{align}
&ta_1=n_1(-a_{10}),\\
&(1-t)a_1=n_2(-a_{10}-a_9-a_3),\\
&t(t-1)a_1=n_3(-a_{10}-t a_9-t^2 a_3).
\end{align}
From the first equation, we obtain $a_{10}=-\frac{t}{n_1}a_1$.
Next, substituting this into the second and the third equations, we obtain
\begin{align}
\begin{split}\label{asdf}
&(1-t)a_1=n_2 \left(\frac{t}{n_1}a_1-a_9-a_3 \right),
\end{split}\\
\begin{split}\label{asdfg}
&(t-1)a_1=n_3 \left(\frac{1}{n_1}a_1-a_9-t a_3 \right).
\end{split}
\end{align}
By calculating $n_3 \times \eqref{asdf}-n_2 \times \eqref{asdfg}$, we obtain
\begin{equation}
(n_2+n_3)(1-t)a_1=\frac{n_2 n_3}{n_1}(t-1)a_1-n_2n_3(1-t)a_3.
\end{equation}
We summarize that
\begin{equation}
(n_1n_2+n_1n_3+n_2n_3)a_1=-n_1n_2n_3a_3.
\end{equation}
Solving on $a_3$, we obtain
\begin{equation}\label{as}
a_3=-\frac{n_1n_2+n_1n_3+n_2n_3}{n_1n_2n_3}a_1.
\end{equation}
Next, in the coordinate system $(z_3,w_3)$ we see that $X=\infty$ is a sigular point. By giving the eigenvalues of the matrix of linear approximation around $X=\infty$, we can obtain
\begin{equation}
n_4(2a_1+a_3)=a_1.
\end{equation}
Substituting \eqref{as} into this equation, we can obtain the relation of the eigenvalues $n_i$:
\begin{equation}
\frac{1}{n_1}+\frac{1}{n_2}+\frac{1}{n_3}+\frac{1}{n_4}=2.
\end{equation}
We have completed the proof of Proposition \ref{Proposition4.1}.

The system \eqref{erty} has the following birational symmetries.
\begin{theorem}
The system \eqref{erty} is invariant under the following transformations\rm{:\rm} with the notation $(*)=(x,y,t;n_1,n_2,n_3,n_4;\alpha_0,\alpha_1,\ldots,\alpha_4),$
\begin{align*}
        s: (*) \rightarrow &(x+\frac{\alpha_2}{y},y,t;n_1,n_2,n_3,n_4;\alpha_0+\alpha_2-n_3 \alpha_2,\\
&\frac{(-n_1 n_2-n_1 n_3-n_2n_3+2n_1n_2n_3)\alpha_1+(-n_1 n_2-n_1 n_3-n_2n_3+n_1n_2n_3)\alpha_2}{-n_1 n_2-n_1 n_3-n_2n_3+2n_1n_2n_3},\\
&-\alpha_2,\alpha_3+\alpha_2-n_2 \alpha_2,\alpha_4+\alpha_2-n_1 \alpha_2), \\
        \pi_1: (*) \rightarrow &(1-x,-y,1-t;n_2,n_1,n_3,n_4;\alpha_0,\alpha_1,\alpha_2,\alpha_4,\alpha_3), \\
        \pi_2: (*) \rightarrow &\left(\frac{t-x}{t-1},-(t-1)y,\frac{t}{t-1};n_3,n_2,n_1,n_4;\alpha_4,\alpha_1,\alpha_2,\alpha_3,\alpha_0 \right), \\
        \pi_3: (*) \rightarrow &(\frac{1}{x},-(yx+\alpha_2)x,\frac{1}{t};n_4,n_2,n_3,n_1;\alpha_0,\frac{\alpha_4 n_2 n_3 n_4}{n_1(2n_1n_2n_3-n_1 n_2-n_1 n_3-n_2n_3)},\\
&\alpha_2,\alpha_3,\frac{\alpha_1 (2n_1n_2n_3-n_1 n_2-n_1 n_3-n_2n_3)}{n_1 n_2 n_3}).
        \end{align*}
\end{theorem}
All transformations satisfy the relation: $s^2={\pi_j}^2=1$. The transformations $\pi_j$ change the eigenvalues $n_1,n_2,n_3,n_4$ in addition to some parameter's changes.

The transformation $s$ is a generalization of the Euler transfomation of the Painlev\'e VI system. The transformations $\pi_j$ correspond to the permutation of the singular points $0,1,t,\infty$. The transformations on sign change of exponents can not be found.

We remark that all transformations coincide with the ones in the case of Painlev\'e VI system when $n_1=n_2=n_3=n_4=2$.

\vspace{0.5cm}
For the system \eqref{erty}, we consider the following problem.
\begin{problem}
When does the system \eqref{erty} have the Painlev\'e property?
\end{problem}
In order to have no movable branch points, the eigenvalues $n_i$ must belong to ${\mathbb Z}$ (see Section 2). At first, let us classify the natural number solutions for the equation \eqref{234}.

\begin{proposition}\label{Proposition4.2}
For the equation \eqref{234} the natural number solutions
$$
\{(n_1,n_2,n_3,n_4) \in {\mathbb N}^4|1 \leq n_1 \leq n_2 \leq n_3 \leq n_4\}
$$
can be classified into four types:
\begin{equation}
\{(n_1,n_2,n_3,n_4)=(1,2,3,6),(1,2,4,4),(1,3,3,3),(2,2,2,2)\}.
\end{equation}
\end{proposition}
We remark that from the symmery of the equation \eqref{234}, we can set 
$$
1 \leq n_1 \leq n_2 \leq n_3 \leq n_4.
$$
The type of $(n_1,n_2,n_3,n_4)=(2,2,2,2)$ is the case of Painlev\'e VI.

{\bf Proof of Proposition \ref{Proposition4.2}.}

(i) The case of $n_1 \geq 3$

By assumption, we see that $\frac{1}{n_i} \leq \frac{1}{3}$. Then, we see that
$$
1\frac{1}{3}=\frac{1}{3}+\frac{1}{3}+\frac{1}{3}+\frac{1}{3} \geq \sum_{i=1}^{4} \frac{1}{n_i}.
$$
This contradicts the equation \eqref{234}.

(ii) The case of $n_1=2$

In this case, we consider
\begin{equation}\label{2345}
\frac{1}{n_2}+\frac{1}{n_3}+\frac{1}{n_4}=\frac{3}{2}.
\end{equation}

(ii-1) The case of $n_2=2$

In this case, we consider
$$
\frac{1}{n_3}+\frac{1}{n_4}=1.
$$
Since $n_j \geq 2$, we see that $\frac{1}{n_j} \leq \frac{1}{2}$. Then, we obtain
$$
n_3=n_4=2.
$$
Consequently, we can obtain $(n_1,n_2,n_3,n_4)=(2,2,2,2)$.

(ii-1) The case of $n_2 \geq 3$

Since $n_j \geq 3$, we see that $\frac{1}{n_j} \leq \frac{1}{3}$. Then, we obtain
$$
\sum_{i=2}^{4} \frac{1}{n_i} \leq 1.
$$
This contradicts the equation \eqref{2345}.

(iii) The case of $n_1=1$

In this case, we consider
\begin{equation}\label{234567}
\frac{1}{n_2}+\frac{1}{n_3}+\frac{1}{n_4}=1.
\end{equation}

(iii-1) The case of $n_2=1$

In this case, we consider
$$
\frac{1}{n_3}+\frac{1}{n_4}=0.
$$
Since $n_j \geq 1$, this is contradiction.

(iii-2) The case of $n_2=2$

In this case, we consider
\begin{equation}\label{23456}
\frac{1}{n_3}+\frac{1}{n_4}=\frac{1}{2}.
\end{equation}
Since $n_3=2$, $n_4 \geq 2$. This is contradiction. So, $n_3 \geq 3$.

If $n_3=3$, we obtain $n_4=6$. Consequently, we can obtain $(n_1,n_2,n_3,n_4)=(1,2,3,6)$.

If $n_3=4$, we obtain $n_4=4$. Consequently, we can obtain $(n_1,n_2,n_3,n_4)=(1,2,4,4)$.

Now, if $5 \leq n_3 \leq n_4$, then $\frac{1}{n_j} \leq \frac{1}{5}$. We obtain
$$
\frac{1}{n_3}+\frac{1}{n_4} \leq \frac{2}{5} < \frac{1}{2}.
$$
This contradicts the equation \eqref{23456}.

(iii-3) The case of $3 \leq n_2 \leq n_3 \leq n_4$

In this case,  from $\frac{1}{n_j} \leq \frac{1}{3}$, we obtain
$$
\frac{1}{n_2}+\frac{1}{n_3}+\frac{1}{n_4} \leq 1.
$$
From the equation \eqref{234567}, we obtain $n_2=n_3=n_4=3$. Consequently, we can obtain $(n_1,n_2,n_3,n_4)=(1,3,3,3)$.

If $4 \leq n_2 \leq n_3 \leq n_4$, then we obtain
$$
\frac{1}{n_2}+\frac{1}{n_3}+\frac{1}{n_4} \leq \frac{3}{4}.
$$
This contradicts the equation \eqref{234567}. Thus, we have completed the proof of Proposition \ref{Proposition4.2}.

\begin{proposition}
For the equation \eqref{234} the integer solutions
$$
(n_1,n_2,n_3,n_4) \in {\mathbb Z}^4 \
$$
can be classified into the following canonical five types\rm{: \rm}
\begin{equation}\label{index1}
\{(n_1,n_2,n_3,n_4)=(1,2,3,6),(1,2,4,4),(1,3,3,3),(2,2,2,2),(n_1,-n_1,1,1)  \ (n_1 \in {\mathbb N}) \}.
\end{equation}
\end{proposition}

\begin{center}
\begin{tabular}{|c||c|} \hline 
Painlev\'e classification & Eq.\eqref{index1} \\ \hline
$k_i=1-\frac{1}{n_i}$ & $n_i$ \\ \hline
Type I:$\left(\frac{m+1}{m},\frac{m-1}{m} \right)$  & $(-m,m,1,1)$ \\ \hline
Type III:$\left(\frac{1}{2},\frac{1}{2},\frac{1}{2},\frac{1}{2} \right)$ & $(2,2,2,2)$ \\ \hline
Type IV:$\left(\frac{2}{3},\frac{2}{3},\frac{2}{3} \right)$ & $(1,3,3,3)$ \\ \hline
Type V:$\left(\frac{3}{4},\frac{3}{4},\frac{1}{2} \right)$ & $(1,4,4,2)$ \\ \hline
Type VI:$\left(\frac{5}{6},\frac{2}{3},\frac{1}{2} \right)$ & $(1,6,3,2)$ \\ \hline
\end{tabular}
\end{center}
In the Painlev\'e classification, Type II is omitted because it may be regarded as a degenerate case of Type III (see \cite{Ince} P 323). Here, we review the Painlev\'e exponent $k$ (see \cite{Ince} P 322);

\begin{align*}
\frac{d}{dz}\begin{pmatrix}
             W \\
             P 
             \end{pmatrix}&=\begin{pmatrix}
             P \\
             \frac{k P^2}{W} 
             \end{pmatrix}=\frac{1}{W} \begin{pmatrix}
             1 & 0  \\
             0 & k
             \end{pmatrix}\begin{pmatrix}
             W P \\
             P^2 
             \end{pmatrix}.
             \end{align*}
If $k \not= 1$, $W(z)=(A z+B)^{\frac{1}{1-k}} \ (A,B \in {\mathbb C})$. In the case of $k=1$, $W(z)=e^{A z+B}$.

Setting $\frac{1}{1-k}=n \ (n \in {\mathbb Z})$, we see the following relation:
\begin{center}
\begin{tabular}{|c||c|} \hline 
Painlev\'e exponent & local index \\ \hline
$k=1-\frac{1}{n}$ & $n$ \\ \hline
\end{tabular}
\end{center}
We remark that $n \in {\mathbb Z}$ for necessary condition of Painlev\'e property. We also note that in \cite{Ince} Page 322 $\frac{1}{1+k}=n \ ( k=1+\frac{1}{n})$ were given.

Let us consider the case $(n_1,n_2,n_3,n_4)=(n_1,-n_1,1,1), \ n_1 \in {\mathbb N}$. In this case, by $\alpha$-method, in order to have no movable branch points at $X=t$ and $X=\infty$, we must put $\alpha_0=\alpha_1=0$ (see relation \eqref{eq:redlocal}) and obtain
\begin{equation}
\begin{pmatrix}
X=0 & X=1 & X=t & X=\infty\\
f_0\begin{pmatrix}
n_1 & \alpha_4\\
0 & 1 
\end{pmatrix} & f_1\begin{pmatrix}
-n_1 & \alpha_3\\
0 & 1 
\end{pmatrix} & f_2\begin{pmatrix}
1 & 0\\
0 & 1 
\end{pmatrix} & f_3\begin{pmatrix}
1 & 0\\
0 & 1 
\end{pmatrix}
\end{pmatrix}.
\end{equation}

The differential system is given by
\begin{equation}\label{system2}
  \left\{
  \begin{aligned}
   \delta t(t-1) \frac{dx}{dt} =&n_1 x(x-1)(x-t)y+n_1 \alpha_1 x^2+\{\alpha_3(t-1)-\alpha_4 t-n_1 \alpha_2 \}x+\alpha_4 t,\\
\delta t(t-1) \frac{dy}{dt} =&-\{ n_1 x^2-(n_1+1)x+t \}y^2+\{-n_1 \alpha_2 x+(n_1(1-t)+1)\alpha_2-\alpha_3(t-1)+\alpha_4 t \}y,
   \end{aligned}
  \right. 
\end{equation}
where $\delta:=n_1 \alpha_2+\alpha_3-\alpha_4$.

Elimination of $y$ from this system gives the second-order ordinary differential equation for the variable $x$; namely,
\begin{equation}\label{equation1}
  \left\{
  \begin{aligned}
   \frac{d^2x}{dt^2} =&\left\{\frac{n_1+1}{n_1(x-1)}+\frac{n_1-1}{n_1 x} \right\} \left(\frac{dx}{dt} \right)^2\\
   &+\frac{\frac{dx}{dt}}{n_1 t(t-1)(x-t)(x-1)x(n_1 \alpha_2+\alpha_3-\alpha_4)} \times\\
   &[n_1(2n_1 t+1-n_1)(t-x)(x-1)x \alpha_2+(t-x)x\{(3t-2)n_1 x+2t(1-n_1)+n_1-2\}\alpha_3\\
   &-(t-x)(x-1)\{(3t-1)n_1 x+(2-n_1)t\}\alpha_4]\\
   &+\frac{1}{n_1(t-1)^2 t^2 (t-x)(x-1)x(n_1 \alpha_2+\alpha_3-\alpha_4)^2} \times\\
   &[(t-1)(t-x)x^2\{n_1 t(x-1)-1+t\}\alpha_3^2-t(t-1)(t-x)(x-1)x(2n_1 x+2-n_1)\alpha_3 \alpha_4\\
   &+t(t-x)(x-1)^2 \{n_1(t-1)x+t \}\alpha_4^2\\
   &+\alpha_2\{n_1(t-1)(n_1 t+1)(t-x)(x-1)x^2 \alpha_3-n_1 t(n_1(t-1)+1)(t-x)(x-1)^2 x \alpha_4 \}].
   \end{aligned}
  \right. 
\end{equation}

Next, we can solve the system \eqref{system2} explicitly. At first, setting
\begin{equation}
  \left\{
  \begin{aligned}
   x_1 =&(x-t)y,\\
y_1 =&\frac{1}{y},
   \end{aligned}
  \right. 
\end{equation}
we can rewrite  the system \eqref{system2}:
\begin{equation}\label{system3}
  \left\{
  \begin{aligned}
   \delta t(t-1) \frac{dx_1}{dt} =&x_1(x_1+\alpha_2),\\
\delta t(t-1) \frac{dy_1}{dt} =&n_1 x_1(x_1+\alpha_2)y_1^2+\{(2n_1 t-1-n_1)x_1+(2n_1 \alpha_2+\alpha_3-\alpha_4)t-(n_1+1)\alpha_2-\alpha_3\}y_1\\
&+n_1 t(t-1),
   \end{aligned}
  \right. 
\end{equation}
where $\delta:=n_1 \alpha_2+\alpha_3-\alpha_4$.

This system is a Riccati extension of the Riccati equation in the variable $x_1$.

\section{The case of Painlev\'e V system}
In this section, we give the Painlev\'e scheme of the Painlev\'e V system. In this case, the accessible singular point $X=Y=0$ has multiplicity of order 2. By making two times blowing-ups, this accessible singular point transformes into simple singular point. For this simple point, we give a matrix of linear approximation around this point. In Proposition \ref{Proposition5.1}, we will show that the condition of the double point  $X=0 \ (2)$ is equivalent to the pair of a simple accessible singular point and a matrix of  degenerate type as linear approximation around this point.
\begin{theorem}\label{Theorem5.1}
For the system \eqref{4}, we give the following Painlev\'e scheme:
\begin{equation}\label{asdjkl}
\begin{pmatrix}
X=0 \ (2) & X=1 & X=\infty\\
$\(\overbrace{\begin{pmatrix}
X'\\
Y' 
\end{pmatrix}=\begin{pmatrix}
0\\
-t
\end{pmatrix}, \ \ f_0\begin{pmatrix}
1 & 0\\
2\alpha_3 & n_3 
\end{pmatrix} }\)$ & f_1\begin{pmatrix}
n_1 & \alpha_0\\
0 & 1 
\end{pmatrix} & f_2\begin{pmatrix}
n_2 & \alpha_1\\
0 & 1 
\end{pmatrix}
\end{pmatrix}.
\end{equation}
Here, $X=0,1,\infty$ are accessible singular points, $f_i \in {\mathbb C}(t)$, $n_i \in {\mathbb C}$, $t \in {\mathbb C}-\{0\}$ and $\alpha_i$ are constant parameters. The symbol $X=0 \ (2)$ means that the point $X=Y=0$ has multiplicity of order 2, and $(X',Y')=(x,x^2 y)$. 
Then, this system coincides with
\begin{equation}\label{erty1}
  \left\{
  \begin{aligned}
   \delta \frac{dx}{dt} =&2n_1 n_2 x^3 y-2n_1 n_2 x^2 y-2n_1(\alpha_1-n_2 \alpha_2)x^2+\{2n_2 \alpha_0+2n_1 \alpha_1-2n_1 n_2 \alpha_2+(n_1+n_2)t\}x\\
&-(n_1+n_2)t,\\
   \delta \frac{dy}{dt} =&-2n_1(2n_2-1)x^2 y^2+2(2n_1 n_2-n_1-n_2)x y^2+2n_1\{2\alpha_1-(3n_2-2)\alpha_2\}x y\\
&-\{2n_2 \alpha_0+2n_1 \alpha_1-2(2n_1 n_2-n_1-n_2)\alpha_2+(n_1+n_2)t\}y+2n_1 \alpha_2(\alpha_1+\alpha_2-n_2 \alpha_2),
   \end{aligned}
  \right. 
\end{equation}
where $\delta:=t\{2n_2 \alpha_0+2n_1 \alpha_1-2(n_1+n_2)\alpha_2+2(2n_1 n_2-n_1-n_2)\alpha_3-(n_1-n_2)t\}$.
\end{theorem}
This system can be considered as a generalization of the Painlev\'e V system. The case of $(n_1,n_2,n_3)=(2,2,2)$ is equivalent to the Painlev\'e V system.

\begin{figure}
\unitlength 0.1in
\begin{picture}( 24.7000, 37.6000)( 25.1000,-40.3000)
%
\special{pn 8}%
\special{pa 2520 530}%
\special{pa 4910 530}%
\special{fp}%
%
\special{pn 20}%
\special{sh 0.600}%
\special{ar 2900 530 22 26  0.0000000 6.2831853}%
%
\special{pn 20}%
\special{pa 2920 1110}%
\special{pa 2920 740}%
\special{fp}%
\special{sh 1}%
\special{pa 2920 740}%
\special{pa 2900 808}%
\special{pa 2920 794}%
\special{pa 2940 808}%
\special{pa 2920 740}%
\special{fp}%
%
\special{pn 8}%
\special{pa 2510 1536}%
\special{pa 4900 1536}%
\special{fp}%
%
\special{pn 20}%
\special{sh 0.600}%
\special{ar 2890 1536 22 26  0.0000000 6.2831853}%
%
\special{pn 20}%
\special{pa 2910 2820}%
\special{pa 2910 2450}%
\special{fp}%
\special{sh 1}%
\special{pa 2910 2450}%
\special{pa 2890 2518}%
\special{pa 2910 2504}%
\special{pa 2930 2518}%
\special{pa 2910 2450}%
\special{fp}%
\put(30.3000,-10.5000){\makebox(0,0)[lb]{Blow up}}%
\put(30.3000,-27.4500){\makebox(0,0)[lb]{Blow up}}%
\put(28.7000,-4.4000){\makebox(0,0)[lb]{$X=Y=0$:double point}}%
\put(30.0000,-14.3000){\makebox(0,0)[lb]{$X_1=Y_1=0$:simple point}}%
%
\special{pn 8}%
\special{pa 2510 3096}%
\special{pa 4900 3096}%
\special{fp}%
%
\special{pn 8}%
\special{pa 2980 1340}%
\special{pa 2590 2130}%
\special{dt 0.045}%
%
\special{pn 8}%
\special{pa 3020 2900}%
\special{pa 2640 3660}%
\special{dt 0.045}%
%
\special{pn 8}%
\special{pa 2610 3490}%
\special{pa 3270 4030}%
\special{fp}%
%
\special{pn 20}%
\special{sh 0.600}%
\special{ar 2810 3340 22 26  0.0000000 6.2831853}%
\put(28.9000,-34.3000){\makebox(0,0)[lb]{$X_3=0,Y_3=-t$:simple point}}%
\put(49.8000,-6.0000){\makebox(0,0)[lb]{$H \cong {\mathbb P}^1$}}%
\end{picture}%
\label{fig:GRS2}
\caption{Resolution of multiplicity of order 2}
\end{figure}
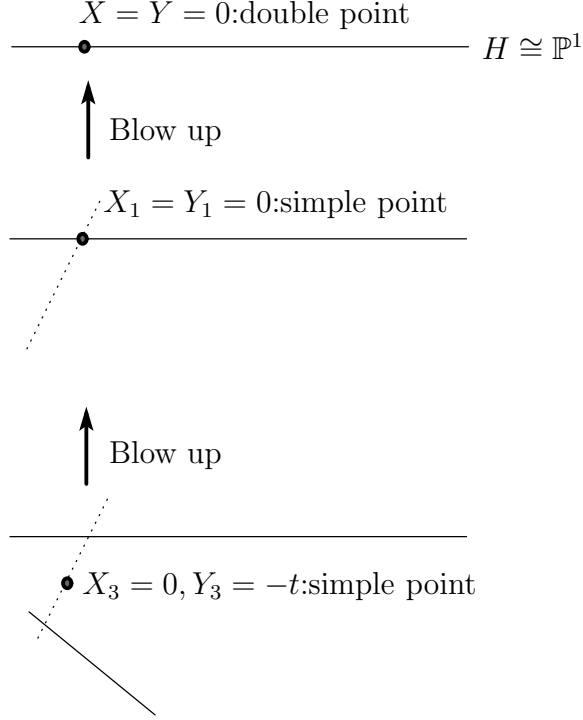

{\bf Proof of Theorem \ref{Theorem5.1}.} \quad We only consider the case of multiplicity of order 2. At first, we can rewrite the system \eqref{4} in the coordinate system $(X,Y)=(x,1/y)$ centered at $(X,Y)=(0,0)$
\begin{equation}
  \left\{
  \begin{aligned}
   \frac{dX}{dt} &=\frac{a_1X^3+a_2X^2+a_5X+a_7}{Y}+\frac{1}{2}\{(3a_1+2a_3)\alpha_2-a_4\}X^2+\{(a_2+a_9)\alpha_2-a_6\}X+a_8,\\
   \frac{dY}{dt} &=-a_{10}-a_9X-a_3X^2-a_4XY-a_6Y-\frac{1}{2}(a_1 \alpha_2+a_4)\alpha_2Y^2 \quad (a_i \in {\mathbb C}(t)).
   \end{aligned}
  \right. 
\end{equation}
By Definition \ref{Def1}, we can calculate the accessible singular points
\begin{equation}
Y=0, \quad a_1X^3+a_2X^2+a_5X+a_7=0.
\end{equation}
By the assumption, $X=Y=0$ is a solution of the system \eqref{6}. Thus, we obtain the condition
$$
a_7=0.
$$
Moreover, since this singular point has multiplicity of order 2, we need the condition
$$
a_5=0.
$$
This condition is necessary condition for multiplicity of order 2.

Next, let us resolve the multiplicity of this point by making two times blowing-ups.

{\bf Step 1.} We blow up at the point $X=Y=0$:
\begin{equation}
X_1=X, \quad Y_1=\frac{Y}{X}.
\end{equation}
Since $X_1=Y_1=0$ must be a singular point, we need the condition
$$
a_{10}=0.
$$
We summarize that the singular point $X=Y=0$ has multiplicity of order 2 if and only if
\begin{equation}\label{sdfg}
a_5=a_7=a_{10}=0.
\end{equation}
{\bf Step 2.} We blow up at the point $X_1=Y_1=0$
\begin{equation}
X_2=X_1, \quad Y_2=\frac{Y_1}{X_1}.
\end{equation}
Here, in order to take a suitable coordinate system we make a change of variables:
\begin{equation}
X_3=X_2, \quad Y_3=\frac{1}{Y_2}.
\end{equation}
We see that the patching data between $(X_3,Y_3)$ and $(x,y)$ is given by $(X_3,Y_3)=(x,x^2 y)$.

In the coordinate system $(X_3,Y_3)$ we rewrite the system \eqref{4}:
\begin{equation}
  \left\{
  \begin{aligned}
   \frac{dX_3}{dt} =&a_8+\frac{1}{2}\{(3a_1+2a_3)\alpha_2-a_4\}X_3^2+a_1X_3Y_3+\{(a_2+a_9)\alpha_2-a_6\}X_3+a_2Y_3,\\
   \frac{dY_3}{dt} =&\frac{2a_8 Y_3}{X_3}+\frac{(2a_2+a_9)Y_3^2}{X_3}+\frac{1}{2}\alpha_2(\alpha_2 a_1+a_4)X_3^2+(2a_1+a_3)Y_3^2+\alpha_2(3a_1+2a_3)X_3 Y_3\\
&+\{2\alpha_2(a_2+a_9)-a_6\}Y_3.
   \end{aligned}
  \right. 
\end{equation}
By the assumption, $(X_3,Y_3)=(0,-t)$ is a simple singular point. So, we obtain  the condition
$$
a_8=\frac{1}{2}t(2a_2+a_9).
$$
Finally, by the assumption, the matrix of linear approximation around $(X_3,Y_3)=(0,-t)$ is given by
\begin{equation}
f_0\begin{pmatrix}
1 & 0\\
2\alpha_3 & n_3 
\end{pmatrix}.
\end{equation}
So, we obtain
$$
a_2=-\frac{1}{4}(n_3+2)a_9, \quad a_1=-\frac{2+2t^2 a_3+t\{2a_6+(n_3-2)\alpha_2 a_9-2\alpha_3 a_9\}}{4t^2}.
$$
For the remaining singular points $X=1$ and $X=\infty$, we can obtain the conditions in the same way of Painlev\'e VI case. Thus, we have completed the proof of Theorem \ref{Theorem5.1}.

By using the conditions \eqref{sdfg} in the proof of Theorem \ref{Theorem5.1}, we easily see the following
\begin{proposition}\label{Proposition5.1}
The condition $X=0 \ (2)$ is equivalent to the following:
\begin{equation}\label{cond1}
\begin{pmatrix}
X=0 \\ 
f \begin{pmatrix}
0 & *\\
0 & 0
\end{pmatrix}
\end{pmatrix},
\end{equation}
\end{proposition}
where $f,* \in {\mathbb C}(t)$.

\begin{remark}
The condition \eqref{cond1} means that $X=0$ is a singular point, and the matrix of linear approximation around this point is given by
\begin{equation}
f \begin{pmatrix}
0 & *\\
0 & 0
\end{pmatrix}.
\end{equation}
\end{remark}

\begin{proposition}\label{Proposition5.2}
The eigenvalues $n_i$ satisfy the following relation:
\begin{equation}\label{555222}
2n_1 n_2 n_3-(n_1+n_2)n_3-2(n_1+n_2)=0.
\end{equation}
\end{proposition}
We see that the case of $(n_1,n_2,n_3)=(2,2,2)$ is equivalent to the Painlev\'e V system.

\begin{proposition}\label{Proposition5.3}
For the equation \eqref{555222} the natural number solutions
$$
\{(n_1,n_2,n_3) \in {\mathbb N}^3|n_1 \geq n_2\}
$$
can be classified into six types:
\begin{equation}
\{(n_1,n_2,n_3)=(2,1,6),(2,2,2),(3,1,4),(3,3,1),(5,1,3),(6,2,1)\}.
\end{equation}
\end{proposition}
We remark that from the symmery of the equation \eqref{555222}, we can set 
$$
n_1 \geq n_2.
$$
The type of $(n_1,n_2,n_3)=(2,2,2)$ is the case of Painlev\'e V.

\vspace{0.2cm}
It is still an open question whether we classify all integer solutions for the equation \eqref{555222}.

\vspace{0.2cm}
Finally, we show that the system \eqref{erty1} has the following birational symmetries.
\begin{theorem}
The system \eqref{erty1} is invariant under the following transformations\rm{:\rm} with the notation $(*)=(x,y,t;n_1,n_2,n_3;\alpha_0,\alpha_1,\alpha_2,\alpha_3),$
\begin{align*}
        s: (*) \rightarrow &(x+\frac{\alpha_2}{y},y,t;n_1,n_2,n_3;\alpha_0+\alpha_2-n_1 \alpha_2,\alpha_1+\alpha_2-n_2 \alpha_2,-\alpha_2,\\
&\frac{2(\alpha_2+\alpha_3)n_1 n_2-(3\alpha_2+\alpha_3)n_1-(3\alpha_2+\alpha_3)n_2}{2n_1n_2-n_1-n_2}),\\
        \pi: (*) \rightarrow &\left(\frac{x}{x-1},-(x-1)((x-1)y+\alpha_2),-t;n_2,n_1,n_3;\alpha_1,\alpha_0,\alpha_2,\alpha_3 \right).
        \end{align*}
\end{theorem}
All transformations satisfy the relation: $s^2={\pi}^2=1$. The transformation $\pi$ changes the eigenvalues $n_1,n_2,n_3$ in addition to some parameter's changes.

The transformation $\pi$ corresponds to the permutation of the singular points $1$ and $\infty$. The transformations on sign change of exponents can not be found.

We remark that all transformations coincide with the ones in the case of Painlev\'e V system when $n_1=n_2=n_3=2$.

\section{The case of Painlev\'e IV system}
In this section, we give the Painlev\'e scheme of the Painlev\'e IV system. In this case, the accessible singular point $X=Y=0$ has multiplicity of order 3. By making three times blowing-ups, this accessible singular point transformes into a simple singular point. For this simple point, we give a matrix of linear approximation around this point. In Proposition \ref{Proposition6.1}, we will show that the condition of the triple point  $X=0 \ (3)$ is equivalent to the pair of a simple accessible singular point and the eigenvalues of two matrices for the expansion (see Proposition \ref{Proposition6.1}) around this point.
\begin{theorem}\label{Theorem6.1}
For the system \eqref{4}, we give the following Painlev\'e scheme:
\begin{equation}\label{asdfgh}
\begin{pmatrix}
X=0 \ (3) & X=\infty\\
$\(\overbrace{\begin{pmatrix}
X'\\
Y' 
\end{pmatrix}=\begin{pmatrix}
0\\
-\frac{1}{2}
\end{pmatrix}, \ \ f_0\begin{pmatrix}
1 & 0\\
2t & n_2 
\end{pmatrix} }\)$ & f_1\begin{pmatrix}
n_1 & \alpha_1\\
0 & 1 
\end{pmatrix}
\end{pmatrix}.
\end{equation}
Here, $X=0,\infty$ are accessible singular points, $f_i \in {\mathbb C}(t)$, $n_i \in {\mathbb C}$, $t \in {\mathbb C}$ and $\alpha_i$ are constant parameters. The symbol $X=0 \ (3)$ means that the point $X=Y=0$ has multiplicity of order 3, and $(X',Y')=(x,x^3 y)$. 
Then, this system coincides with
\begin{equation}\label{erty2}
  \left\{
  \begin{aligned}
   \frac{dx}{dt} =&a(t)\left(x^3 y+\frac{(n_1 \alpha_2-\alpha_1)x^2}{n_1}+\frac{(2n_1-1)t x}{3n_1}+\frac{n_1+1}{6n_1} \right),\\
   \frac{dy}{dt} =&a(t)\left(-\frac{(2n_1-1)x^2 y^2}{n_1}+\frac{(2\alpha_1-(3n_1-2)\alpha_2)}{n_1}xy-\frac{(2n_1-1)ty}{3n_1}+\frac{\alpha_2(\alpha_1-(n_1-1)\alpha_2)}{n_1} \right),
   \end{aligned}
  \right. 
\end{equation}
where $a(t) \in {\mathbb C}(t)$.
\end{theorem}
This system can be considered as a generalization of the Painlev\'e IV system. The case of $(n_1,n_2)=(2,3)$ and $a(t)=4$ is equivalent to the Painlev\'e IV system.

\begin{figure}
\unitlength 0.1in
\begin{picture}( 62.1000, 28.1000)( 11.4000,-30.3000)
%
\special{pn 8}%
\special{pa 1290 490}%
\special{pa 3680 490}%
\special{fp}%
%
\special{pn 20}%
\special{sh 0.600}%
\special{ar 1670 490 22 26  0.0000000 6.2831853}%
%
\special{pn 8}%
\special{pa 4960 496}%
\special{pa 7350 496}%
\special{fp}%
%
\special{pn 20}%
\special{sh 0.600}%
\special{ar 5340 496 22 26  0.0000000 6.2831853}%
\put(41.8000,-10.2500){\makebox(0,0)[lb]{Blow up}}%
\put(16.4000,-4.0000){\makebox(0,0)[lb]{$X=Y=0$:triple point}}%
\put(54.5000,-3.9000){\makebox(0,0)[lb]{$X_1=Y_1=0$:double point}}%
%
\special{pn 8}%
\special{pa 1880 1706}%
\special{pa 4270 1706}%
\special{fp}%
%
\special{pn 8}%
\special{pa 5430 300}%
\special{pa 5040 1090}%
\special{dt 0.045}%
%
\special{pn 8}%
\special{pa 2390 1510}%
\special{pa 2010 2270}%
\special{dt 0.045}%
%
\special{pn 8}%
\special{pa 1980 2100}%
\special{pa 2640 2640}%
\special{fp}%
%
\special{pn 20}%
\special{sh 0.600}%
\special{ar 2300 1700 22 26  0.0000000 6.2831853}%
\put(24.0000,-15.9000){\makebox(0,0)[lb]{$X_2=Y_2=0$:simple point}}%
\put(37.5000,-5.6000){\makebox(0,0)[lb]{$H \cong {\mathbb P}^1$}}%
%
\special{pn 8}%
\special{pa 4940 1766}%
\special{pa 7330 1766}%
\special{fp}%
%
\special{pn 8}%
\special{pa 5450 1570}%
\special{pa 5070 2330}%
\special{dt 0.045}%
%
\special{pn 8}%
\special{pa 5040 2160}%
\special{pa 5700 2700}%
\special{fp}%
%
\special{pn 20}%
\special{sh 0.600}%
\special{ar 5240 2020 22 26  0.0000000 6.2831853}%
\put(53.8000,-21.2000){\makebox(0,0)[lb]{$X_4=0,Y_4=-\frac{1}{2}$:simple point}}%
%
\special{pn 8}%
\special{pa 5620 2420}%
\special{pa 5180 3030}%
\special{fp}%
%
\special{pn 20}%
\special{pa 4720 776}%
\special{pa 4090 776}%
\special{fp}%
\special{sh 1}%
\special{pa 4090 776}%
\special{pa 4158 796}%
\special{pa 4144 776}%
\special{pa 4158 756}%
\special{pa 4090 776}%
\special{fp}%
\put(41.6000,-24.5000){\makebox(0,0)[lb]{Blow up}}%
%
\special{pn 20}%
\special{pa 4700 2200}%
\special{pa 4070 2200}%
\special{fp}%
\special{sh 1}%
\special{pa 4070 2200}%
\special{pa 4138 2220}%
\special{pa 4124 2200}%
\special{pa 4138 2180}%
\special{pa 4070 2200}%
\special{fp}%
\put(12.3000,-24.3500){\makebox(0,0)[lb]{Blow up}}%
%
\special{pn 20}%
\special{pa 1770 2186}%
\special{pa 1140 2186}%
\special{fp}%
\special{sh 1}%
\special{pa 1140 2186}%
\special{pa 1208 2206}%
\special{pa 1194 2186}%
\special{pa 1208 2166}%
\special{pa 1140 2186}%
\special{fp}%
\end{picture}%
\label{fig:GRS3}
\caption{Resolution of multiplicity of order 3}
\end{figure}
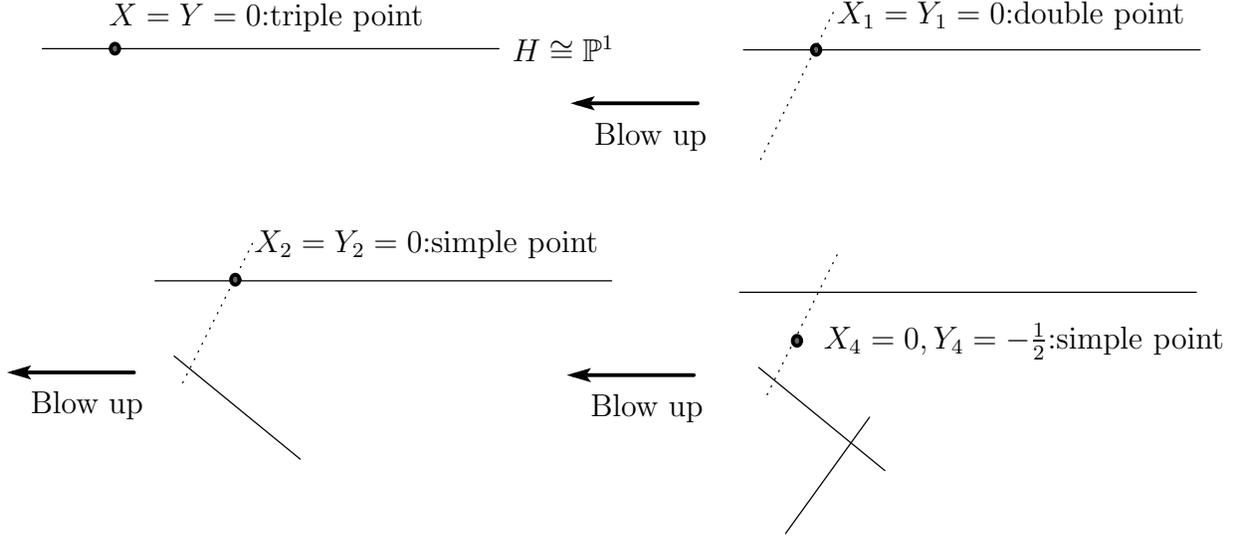

{\bf Proof of Theorem \ref{Theorem6.1}.} \quad We only consider the case of multiplicity of order 3. At first, we can rewrite the system \eqref{4} in the coordinate system $(X,Y)=(x,1/y)$ centered at $(X,Y)=(0,0)$
\begin{equation}\label{putry}
  \left\{
  \begin{aligned}
   \frac{dX}{dt} &=\frac{a_1X^3+a_2X^2+a_5X+a_7}{Y}+\frac{1}{2}\{(3a_1+2a_3)\alpha_2-a_4\}X^2+\{(a_2+a_9)\alpha_2-a_6\}X+a_8,\\
   \frac{dY}{dt} &=-a_{10}-a_9X-a_3X^2-a_4XY-a_6Y-\frac{1}{2}(a_1 \alpha_2+a_4)\alpha_2Y^2 \quad (a_i \in {\mathbb C}(t)).
   \end{aligned}
  \right. 
\end{equation}
By Definition \ref{Def1}, we can calculate the accessible singular points
\begin{equation}
Y=0, \quad a_1X^3+a_2X^2+a_5X+a_7=0.
\end{equation}
By the assumption, $X=Y=0$ is a solution of the system \eqref{6}. Thus, we obtain the condition
$$
a_7=0.
$$
Moreover, since this singular point has multiplicity of order 3, we need the conditions
$$
a_2=a_5=0.
$$
This condition is necessary condition for multiplicity of order 3.

Next, let us resolve the multiplicity of this point by making three times blowing-ups.

{\bf Step 1.} We blow up at the point $X=Y=0$:
\begin{equation}
X_1=X, \quad Y_1=\frac{Y}{X}.
\end{equation}
Since $X_1=Y_1=0$ must be a double singular point, we need the condition
$$
a_{10}=0.
$$
{\bf Step 2.} We blow up at the point $X_1=Y_1=0$
\begin{equation}
X_2=X_1, \quad Y_2=\frac{Y_1}{X_1}.
\end{equation}
Since $X_2=Y_2=0$ must be a singular point, we need the condition
$$
a_{9}=0.
$$
We summarize that the singular point $X=Y=0$ has multiplicity of order 3 if and only if
\begin{equation}\label{sdfggggg}
a_2=a_5=a_7=a_{9}=a_{10}=0.
\end{equation}
{\bf Step 3.} We blow up at the point $X_2=Y_2=0$
\begin{equation}
X_3=X_2, \quad Y_3=\frac{Y_2}{X_2}.
\end{equation}
Here, in order to take a suitable coordinate system we make a change of variables:
\begin{equation}
X_4=X_3, \quad Y_4=\frac{1}{Y_3}.
\end{equation}
We see that the patching data between $(X_4,Y_4)$ and $(x,y)$ is given by $(X_4,Y_4)=(x,x^3 y)$.

By doing the same argument in the proof of Theorem \ref{Theorem5.1}, we can obtain some conditions.

For the remaining singular point $X=\infty$, we can obtain the conditions in the same way of Painlev\'e VI case. Thus, we have completed the proof of Theorem \ref{Theorem6.1}. \qed

By using the conditions \eqref{sdfggggg} in the proof of Theorem \ref{Theorem6.1}, we easily see the following
\begin{proposition}\label{Proposition6.1}
The condition $X=0 \ (3)$ is equivalent to the following conditions:
$X=0$ is a singular point, and the eigenvalues $a_2,a_5,a_9,a_{10}$ for the following expansion of the system \eqref{putry} are given by
$$
a_2=a_5=a_{9}=a_{10}=0,
$$
where
\begin{align}\label{6767654}
\frac{d}{dt}\begin{pmatrix}
             X \\
             Y 
             \end{pmatrix}&=\frac{1}{Y}\left\{\begin{pmatrix}
             a_5 & *  \\
             0 & -a_{10} 
             \end{pmatrix}\begin{pmatrix}
             X \\
             Y 
             \end{pmatrix}+\begin{pmatrix}
             a_2 & *  \\
             0 & -a_{9}
             \end{pmatrix}\begin{pmatrix}
             X^2 \\
             XY 
             \end{pmatrix}+\begin{pmatrix}
             a_1 & *  \\
             0 & -a_3
             \end{pmatrix}\begin{pmatrix}
             X^3 \\
             X^2 Y 
             \end{pmatrix}+\cdots\right\},
             \end{align}
             where  $* \in {\mathbb C}(t)$.
\end{proposition}

\begin{proposition}\label{Proposition6.2}
The eigenvalues $n_i$ satisfy the following relation:
\begin{equation}\label{555222345}
2n_1 n_2-3n_1-n_2-3=0.
\end{equation}
\end{proposition}
We see that the case of $(n_1,n_2)=(2,3)$ is equivalent to the Painlev\'e IV system.

\begin{proposition}\label{Proposition6.3}
For the equation \eqref{555222345} the natural number solutions can be classified into three types:
\begin{equation}
\{(n_1,n_2)=(1,6),(2,3),(5,2)\}.
\end{equation}
\end{proposition}
The type of $(n_1,n_2)=(2,3)$ is the case of Painlev\'e IV.

{\bf Proof of Proposition \ref{Proposition6.3}.} \quad At first, we rewrite the equation \eqref{555222345} as follows:
$$
n_1=\frac{n_2+3}{2n_2-3}.
$$

(i) The case of $n_1 \geq 6$

By assumption, we see that $6 \leq \frac{n_2+3}{2n_2-3}$. Then, we see that
$$
n_2 \leq \frac{21}{11}<1.
$$
This contradicts the condition $n_2 \geq 1$. Thus, we see that $n_1 \leq 5$.

The remaining cases can be solved by the same argument in the proof of Proposition \ref{Proposition4.2}. \qed

\vspace{0.2cm}
It is still an open question whether we classify all integer solutions for the equation \eqref{555222345}.

\vspace{0.2cm}
Finally, we show that the system \eqref{erty2} has the following birational symmetry.
\begin{theorem}
The system \eqref{erty2} is invariant under the following transformation\rm{:\rm} with the notation $(*)=(x,y,t;n_1,n_2;\alpha_0,\alpha_1,\alpha_2),$
\begin{align*}
        s: (*) \rightarrow &\left(x+\frac{\alpha_2}{y},y,t;n_1,n_2;\alpha_1+\alpha_2-n_1 \alpha_2,-\alpha_2 \right).
        \end{align*}
\end{theorem}
The transformation satisfies the relation: $s^2=1$. The transformations on sign change of exponents can not be found.

We remark that the transformation $s$ coincides with the one in the case of Painlev\'e IV system when $(n_1,n_2)=(2,3)$.

\section{The case of Painlev\'e III system}
In this section, we give the Painlev\'e scheme of the Painlev\'e III system. In this case, two accessible singular points have multiplicity of order 2.
\begin{theorem}\label{Theorem7.1}
For the system \eqref{4}, we give the following Painlev\'e scheme:
\begin{equation}\label{asdfghhhhh}
\begin{pmatrix}
X=0 \ (2) & X=\infty \ (2)\\
$\(\overbrace{\begin{pmatrix}
X'\\
Y' 
\end{pmatrix}=\begin{pmatrix}
0\\
-t
\end{pmatrix}, \ \ f_0\begin{pmatrix}
1 & 0\\
2\alpha_0 & n_1
\end{pmatrix} }\)$ & $\(\overbrace{\begin{pmatrix}
X''\\
Y'' 
\end{pmatrix}=\begin{pmatrix}
0\\
-1
\end{pmatrix}, \ \ f_1\begin{pmatrix}
1 & 0\\
2\alpha_1 & n_2 
\end{pmatrix} }\)$
\end{pmatrix}.
\end{equation}
Here, $X=0,\infty$ are accessible singular points, $f_i \in {\mathbb C}(t)$, $n_i \in {\mathbb C}$, $t \in {\mathbb C}-\{0\}$ and $\alpha_i$ are constant parameters, and $(X',Y')=(x,x^2 y)$ and $(X'',Y'')=\left(\frac{1}{x},-\frac{(xy+\alpha_2)}{x} \right)$. 
Then, this system coincides with
\begin{equation}\label{erty3}
  \left\{
  \begin{aligned}
   \delta t \frac{dx}{dt} =&-(n_1+2)x^2 y+2x^2+2(n_1 \alpha_1-2\alpha_2)x-n_1 t,\\
   \delta t \frac{dy}{dt} =&4x y^2-4xy-(2n_1 \alpha_1+(n_1-6)\alpha_2)y-2\alpha_2,
   \end{aligned}
  \right. 
\end{equation}
where $\delta:=4\alpha_0+2n_1 \alpha_1-(n_1+2)\alpha_2$.
\end{theorem}
This system can be considered as a generalization of the Painlev\'e III system. The case of $(n_1,n_2)=(2,2)$ is equivalent to the Painlev\'e III system.

By the same way of Painlev\'e V case, we can prove Theorem \ref{Theorem7.1}. \qed

\begin{proposition}\label{Proposition7.1}
The eigenvalues $n_i$ satisfy the following relation:
\begin{equation}\label{555222345456}
n_1 n_2=4.
\end{equation}
\end{proposition}
We see that the case of $(n_1,n_2)=(2,2)$ is equivalent to the Painlev\'e III system.

\begin{proposition}\label{Proposition7.2}
For the equation \eqref{555222345456} the natural number solutions can be classified into two types:
\begin{equation}
\{(n_1,n_2)=(2,2),(4,1)\}.
\end{equation}
\end{proposition}
The type of $(n_1,n_2)=(2,2)$ is the case of Painlev\'e III.

\vspace{0.2cm}
It is still an open question whether we classify all integer solutions for the equation \eqref{555222345456}.

\vspace{0.2cm}
Finally, we show that the system \eqref{erty3} has the following birational symmetries.
\begin{theorem}
The system \eqref{erty3} is invariant under the following transformations\rm{:\rm} with the notation $(*)=(x,y,t;n_1,n_2;\alpha_0,\alpha_1,\alpha_2),$
\begin{align*}
        s: (*) \rightarrow &\left(x+\frac{\alpha_2}{y},y,t;n_1,n_2;\alpha_0+\alpha_2-n_1 \alpha_2,\alpha_1+\alpha_2-n_2 \alpha_2,-\alpha_2 \right),\\
        \pi: (*) \rightarrow &\left(\frac{t}{x},-\frac{x(xy+\alpha_2)}{t},t;n_2,n_1;\alpha_1,\alpha_0,\alpha_2 \right).
        \end{align*}
\end{theorem}
All transformations satisfy the relation: $s^2={\pi}^2=1$. The transformation $\pi$ changes the eigenvalues $n_1,n_2$ in addition to some parameter's changes.

The transformation $\pi$ corresponds to the permutation of the singular points $0$ and $\infty$. The transformations on sign change of exponents can not be found.

We remark that all transformations coincide with the ones in the case of Painlev\'e III system when $n_1=n_2=2$.

\section{Existence theorem of non-linear ordinary differential systems in dimension two with only simple accessible singular points}

For a linear differential equation of Fuchs type, it is well-known that
\begin{theorem}\label{Exist:th}
Let us consider the n-th order linear ordinary differential equations$:$
\begin{align}\label{lineq}
\begin{split}
\frac{d^nx}{dt^n}+a_1(t)\frac{d^{n-1}x}{dt^{n-1}}+\cdots+a_{n-1}(t)\frac{dx}{dt}+a_n(t) x=0,
\end{split}
\end{align}
where $a_i(t)$ are meromorphic functions defined in a domain in the Riemann sphere ${\mathbb P}^1$.

\noindent
There exists an ordinary differential equation with n-th order satisfying the assumptions $(F1),(F2)$ and $(F3)$.

$(F1)$ This equation has only $(m+1)$ points $x=c_j$ on the Riemann sphere ${\mathbb P}^1$ as its regular singular points.

$(F2)$ Its local exponent at each singular point $c_j$ coincides with ${\rho}_{jl} \in {\mathbb C} \quad (j=1,2,\ldots,m+1,l=1,2,\ldots,n)$.

$(F3)$ ${\rho}_{jl} \ (j=1,2,\ldots,m+1,l=1,2,\ldots,n)$ given in $(F2)$ satisfies the Fuchs' relation$:$
\begin{equation}
\Sigma_{j=1}^{m+1} \Sigma_{l=1}^{n} {\rho}_{jl}=\frac{(m-1)n(n-1)}{2}.
\end{equation}
\end{theorem}

Let us consider the following problem.
\begin{problem}
Can we construct a non-linear ordinary differential system in dimension two satisfying similar conditions of $(F1),(F2)$ and $(F3)$  from the viewpoint of geometrical property? 
\end{problem}

In this section, let us consider a system of the first-order ordinary differential equations in dimension two:
\begin{equation}\label{fq}
  \left\{
  \begin{aligned}
   \frac{dx}{dt} &=f_1(x,y),\\
   \frac{dy}{dt} &=f_2(x,y) \quad (f_i \in {\mathbb C}(t)[x,y]).
   \end{aligned}
  \right. 
\end{equation}
We assume that the regular vector field associated with the system \eqref{fq} defined on ${\mathbb C}^2 \times B$
\begin{equation}
v=\frac{\partial}{\partial t}+\frac{dx}{dt}\frac{\partial}{\partial x}+\frac{dy}{dt}\frac{\partial}{\partial y}
\end{equation}
is extended to to a rational vector field $\tilde v$ on ${\Sigma_n} \times B$
\begin{equation}\label{logcond}
\tilde v \in H^0({\Sigma_n \times B},\Theta_{\Sigma_n \times B}(-\log{D^{(0)}})(D^{(0)})),
\end{equation}
where $B$ is a domain in ${\mathbb C}$.

Here, we review the algebraic surface ${\Sigma_n}$, which is obtained by gluing four copies of ${\mathbb C}^2$ via the following identification.
\begin{align}
\begin{split}
&U_j \cong {\mathbb C}^2 \ni (z_j,w_j) \ (j=0,1,2,3)\\
&z_0=x, \ w_0=y, \quad z_1=\frac{1}{x}, \ w_1=-x^n y-\alpha x,\\
&z_2=z_0, \ w_2=\frac{1}{w_0}, \quad z_3=z_1, \ w_3=\frac{1}{w_1},
\end{split}
\end{align}
where $\alpha$ is a complex constant parameter.

We define a divisor $D^{(0)}$ on ${\Sigma_n}$:
\begin{equation}
D^{(0)}=\{(z_2,w_2) \in U_2|w_2=0\} \cup \{(z_3,w_3) \in U_3|w_3=0\} \cong {\mathbb P}^1.
\end{equation}
The self-intersection number of $D^{(0)}$ is given by
\begin{equation}
(D^{(0)})^2=n.
\end{equation}

The condition \eqref{logcond} is equivalent to the following:
\begin{enumerate}
\item Holomorphy in the coordinate system $(z_1,w_1)=(1/x,-x^n y-\alpha x)$,
\item In the coordinate system $(X,Y)=(x,1/y)$, the differential system \eqref{fq} must be taken of the form:
\begin{equation}
  \left\{
  \begin{aligned}
   \frac{dX}{dt} &=\frac{F_1(X,Y)}{Y},\\
   \frac{dY}{dt} &=F_2(X,Y) \quad (F_i \in {\mathbb C}(t)[X,Y]).
   \end{aligned}
  \right. 
\end{equation}
\end{enumerate}

In the coordinate system $(z_1,w_1)$ the right hand side of this system is polynomial with respect to $z_1,w_1$. However, on the boundary divisor $D^{(0)} \cong {\mathbb P}^1$ this system has a pole in each coordinate system $(z_i,w_i) \ i=2,3$. By rewriting the system at each singular point, this rational vector field has a pole along the divisor $D^{(0)}$, whose order is one.

In this section, we consider the case of {\it simple} accessible singular points.

The following theorem can be considered as a non-linear version of Theorem \ref{Exist:th} from the viewpoint of geometrical property.

\begin{center}
\begin{tabular}{|c||c|c|} \hline 
 & Linear  & Non-linear    \\ \hline
Category & Equation \eqref{lineq} &  $H^0({\Sigma_n \times B},\Theta_{\Sigma_n \times B}(-\log{D^{(0)}})(D^{(0)}))$  \\ \hline 
Condition 1 & regular type & simple accessible type (see Section 2)   \\ \hline 
Condition 2 & local exponent & local index (see Section 2) \\ \hline 
Condition 3 & Fuchs' relation & relation \eqref{relation} \\ \hline
\end{tabular}
\end{center}

\begin{theorem}\label{nonlinear:th}
Let us consider an ordinary differential system in dimension two satisfying the condition \eqref{logcond}. There exists an ordinary differential system of this type satisfying the assumptions $(A1),(A2)$ and $(A3)$.

$(A1)$ This system has only $(n+2)$ points $c_1,c_2,\ldots,c_{n},t,\infty$ on the boundary divisor $D^{(0)} \times B$ as its simple accessible singular points, where $c_i \in {\mathbb C}$ and $t \in B$.

$(A2)$ The ratio of its local index at each accessible singular point $c_i$ coincides with $m_i \in {\mathbb C}-\{0\}$.

$(A3)$ $m_i \ (i=1,2,\ldots,n+2)$ given in $(A2)$ satisfies the relation$:$
\begin{equation}\label{relation}
\Sigma_{i=1}^{n+2} \frac{1}{m_i}=n.
\end{equation}
\end{theorem}
We note that the {\it simple} accessible singular point $P$ means that $P$ is an accessible singular point and has its multiplicity of order 1.

In this paper, we find
\begin{equation}\label{eq:1}
  \left\{
  \begin{aligned}
   \frac{dx}{dt} =&a_1(t)(x-c_1)(x-c_2)\ldots(x-c_n)(x-t)y+b_1[x],\\
   \frac{dy}{dt} =&-\frac{a_1(t)}{m_1m_2 \ldots m_{n+1}}y^2\{m_2m_3 \ldots m_{n+1}(x-c_2)(x-c_3) \ldots (x-c_{n})(x-t)\\
   &+m_1m_3 \ldots m_{n+1}(x-c_1)(x-c_3) \ldots (x-c_{n})(x-t)\\
   &+ \ldots\\
   &+m_1m_2 \ldots m_{n}(x-c_1)(x-c_2) \ldots (x-c_{n-1})(x-c_{n}) \}+b_2[x]y+b_3[x],
   \end{aligned}
  \right. 
\end{equation}
where $a_1(t) \in {\mathbb C}(t)$ and $b_i[x] \in {\mathbb C}(t)[x]$ satisfy certain conditions in order to become a polynomial class in the coordinate system $(z_1,w_1)$.

\begin{center}
\begin{tabular}{|c||c|} \hline
Equation & System \eqref{eq:1}  \\ \hline
Compactification & ${\Sigma}_n$: Hirzebruch surface of degree n  \\ \hline
Accessible singular points & $(z_2,w_2)=(c_i,0) \ (i=1,2,\ldots,n),(t,0),(\infty,0)$  \\ \hline
Painlev\'e scheme & $
\begin{pmatrix}
P_i:z_2=c_i, &P_{n+1}:z_2=t, & P_{n+2}:z_2=\infty\\
\begin{pmatrix}
m_i & * \\
0 & 1 
\end{pmatrix} & \begin{pmatrix}
m_{n+1} & * \\
0 & 1 
\end{pmatrix} & \begin{pmatrix}
m_{n+2} & * \\
0 & 1 
\end{pmatrix}
\end{pmatrix}
$  \\ \hline
Relation of eigenvalues $m_i$ & $\Sigma_{i=1}^{n+2} \frac{1}{m_i}=n$  \\ \hline
\end{tabular}
\end{center}

Before we will show that this system satisfies the assumptions $(A1),(A2)$ and $(A3)$, we can check the conditions of vector field.

(i) Degree of polynomials $f_i(x,y)$ with respect to $y$

If the system \eqref{fq} belongs in $H^0({\Sigma_n \times B},\Theta_{\Sigma_n \times B}(-\log{D^{(0)}})({D^{(0)}}))$,
in the coordinate system $(X,Y)=(x,1/y)$ this system must be taken of the form:
\begin{equation}
  \left\{
  \begin{aligned}
   \frac{dX}{dt} &=\frac{F_1(X,Y)}{Y},\\
   \frac{dY}{dt} &=F_2(X,Y) \quad (F_i \in {\mathbb C}(t)[X,Y]).
   \end{aligned}
  \right. 
\end{equation}
By this condition, we see that the system \eqref{fq} must be taken of the form:
\begin{equation}\label{11111111}
  \left\{
  \begin{aligned}
   \frac{dx}{dt} &=b_1(x)+b_2(x)y,\\
   \frac{dy}{dt} &=b_3(x)+b_4(x)y+b_5(x)y^2 \quad (b_i \in {\mathbb C}(t)[x]).
   \end{aligned}
  \right. 
\end{equation}
Here, the degree of each $b_i$  with respect to $x$ is given by
\begin{equation}
deg(b_1)=l, \ deg(b_2)=m, \ deg(b_3)=p, \ deg(b_4)=q, \ deg(b_5)=r,
\end{equation}
where $l,m,n,p,r \in {\mathbb N}$.

(ii) Holomorphy in the coordinate system

In the coordinate system:
\begin{equation}
(x_1,y_1)=(1/x,-y x^n-g_{n-1} x^{n-1}-\cdots-g_1 x),
\end{equation}
the first equation of the system \eqref{11111111} is given by
\begin{equation}
\frac{d x_1}{dt}=-x_1^2\left\{b_1\left(\frac{1}{x_1} \right)+b_2\left(\frac{1}{x_1} \right)(-y_1 x_1^n-g_1 x_1^{n-1}-\cdots-g_{n-1} x_1) \right\}.
\end{equation}
Since the right hand side of this system must be polynomial with respect to $x_1$, we compare two terms
\begin{equation}
  \left\{
  \begin{aligned}
   b_1\left(\frac{1}{x_1} \right) &=\frac{b_1^{(l)}}{x_1^l}+\cdots,\\
   b_2\left(\frac{1}{x_1} \right) &=-g_{n-1} \frac{b_2^{(m)}}{x_1^{m-1}}+\cdots.
   \end{aligned}
  \right. 
\end{equation}
Since $b_1^{(l)} \not=0$ and $b_2^{(m)} \not=0$, we can obtain
\begin{equation}
l=m-1
\end{equation}
Next, we compare the term involving $y_1$:
\begin{equation}
x_1^{n+2} b_2\left(\frac{1}{x_1} \right)y_1=x_1^{n+2} \left(\frac{b_2^{(m)}}{x_1^m}+\frac{b_2^{(m-1)}}{x_1^{m-1}}+\cdots \right)y_1 \quad (b_2^{(j)} \in {\mathbb C}(t)).
\end{equation}
If this becomes polynomial with respect to $x_1,y_1$, 
\begin{equation}
m=n+2.
\end{equation}
In the same way, we can obtain
\begin{equation}
deg(b_1)=n+1, \ deg(b_2)=n+2, \ deg(b_3)=n-1, \ deg(b_4)=n, \ deg(b_5)=n+1.
\end{equation}

At first, we remark that
\begin{proposition}
These systems \eqref{erty} and \eqref{PVI} satisfy the assumptions $(A1),(A2)$ and $(A3)$.
\end{proposition}

Next, in general case we consider
\begin{equation}\label{fc1}
  \left\{
  \begin{aligned}
   \frac{dx}{dt} =&a_1(t)(x-c_1)(x-c_2)\ldots(x-c_n)(x-t)y+b_1[x],\\
   \frac{dy}{dt} =&-\frac{a_1(t)}{m_1m_2 \ldots m_{n+1}}y^2\{m_2m_3 \ldots m_{n+1}(x-c_2)(x-c_3) \ldots (x-c_{n})(x-t)\\
   &+m_1m_3 \ldots m_{n+1}(x-c_1)(x-c_3) \ldots (x-c_{n})(x-t)\\
   &+ \ldots\\
   &+m_1m_2 \ldots m_{n}(x-c_1)(x-c_2) \ldots (x-c_{n-1})(x-c_{n}) \}+b_2[x]y+b_3[x],
   \end{aligned}
  \right. 
\end{equation}
where $a_1(t) \in {\mathbb C}(t)$ and $b_i[x] \in {\mathbb C}(t)[x]$ satisfy certain conditions in order to become a polynomial class in the coordinate system $(z_1,w_1)$.

In the coordinate system $(z_1,w_1)$, the system \eqref{fc1} can be rewritten as follows:
\begin{equation}\label{fc2}
  \left\{
  \begin{aligned}
   \frac{dz_1}{dt} =&-a_1(t)z_1(1-c_1 z_1)(1-c_2 z_1)\ldots(1-c_n z_1)(1-t z_1)w_1+c_1[z_1],\\
   \frac{dw_1}{dt} =&n a_1(t)(1-c_1 z_1)(1-c_2 z_1)\ldots(1-c_n z_1)(1-t z_1)w_1^2\\
&-\frac{a_1(t)}{m_1m_2 \ldots m_{n+1}}{w_1}^2\{m_2m_3 \ldots m_{n+1}(1-c_2 z_1)(1-c_3 z_1) \ldots (1-c_{n} z_1)(1-t z_1)\\
   &+m_1m_3 \ldots m_{n+1}(1-c_1 z_1)(1-c_3 z_1) \ldots (1-c_{n} z_1)(1-t z_1)\\
   &+ \ldots\\
   &+m_1m_2 \ldots m_{n}(1-c_1 z_1)(1-c_2 z_1) \ldots (1-c_{n-1} z_1)(1-c_{n} z_1) \}+c_2[z_1]w_1+c_3[z_1],
   \end{aligned}
  \right. 
\end{equation}
where $c_i[z_1] \in {\mathbb C}(t)[z_1]$ satisfy certain conditions in order to become a polynomial class in the coordinate system $(x,y)=(1/z_1,-z_1^n w_1-\alpha z_1^{n-1})$.

\begin{proposition}\label{prop1}
The system \eqref{fc1} satisfies the assumptions $(A1),(A2)$ and $(A3)$.
\end{proposition}

{\bf Proof.} In the coordinate system $(X,Y)=(x,1/y)$ the system \eqref{fc1} can be rewritten as follows:
\begin{equation}\label{fc11}
  \left\{
  \begin{aligned}
   \frac{dX}{dt} =&\frac{a_1(t)(X-c_1)(X-c_2)\ldots(X-c_n)(X-t)}{Y}+b_1[X],\\
   \frac{dY}{dt} =&\frac{a_1(t)}{m_1m_2 \ldots m_{n+1}}\{m_2m_3 \ldots m_{n+1}(X-c_2)(X-c_3) \ldots (X-c_{n})(X-t)\\
   &+m_1m_3 \ldots m_{n+1}(X-c_1)(X-c_3) \ldots (X-c_{n})(X-t)\\
   &+ \ldots\\
   &+m_1m_2 \ldots m_{n}(X-c_1)(X-c_2) \ldots (X-c_{n-1})(X-c_{n}) \}-b_2[X]Y-b_3[X]Y^2.
   \end{aligned}
  \right. 
\end{equation}
By Definition \ref{Def1}, we can calculate its accessible singular points
\begin{equation}
Y=0, \quad (X-c_1)(X-c_2)\ldots(X-c_n)(X-t)=0.
\end{equation}
We obtain
\begin{equation}
X=c_1,c_1,\ldots,c_n,t.
\end{equation}

Next, let us calculate its local index at each point. At first, in the coordinate system $(X_1,Y_1)=(X-c_1,Y)$ the system \eqref{fc11} can be rewritten as follows:
\begin{equation}
  \left\{
  \begin{aligned}
   \frac{dX_1}{dt} =&\frac{a_1(t)X_1(X_1+c_1-c_2)\ldots(X_1+c_1-c_n)(X_1+c_1-t)}{Y_1}+b_1[X_1+c_1],\\
   \frac{dY_1}{dt} =&\frac{a_1(t)}{m_1m_2 \ldots m_{n+1}} \times\\
&\{m_2m_3 \ldots m_{n+1}(X_1+c_1-c_2)(X_1+c_1-c_3) \ldots (X_1+c_1-c_{n})(X_1+c_1-t)\\
   &+m_1m_3 \ldots m_{n+1}X_1(X_1+c_1-c_3) \ldots (X_1+c_1-c_{n})(X_1+c_1-t)\\
   &+ \ldots\\
   &+m_1m_2 \ldots m_{n}X_1(X_1+c_1-c_2) \ldots (X_1+c_1-c_{n-1})(X_1+c_1-c_{n}) \}\\
&-b_2[X_1+c_1]Y_1-b_3[X_1+c_1]Y_1^2.
   \end{aligned}
  \right. 
\end{equation}
The matrix of linear approximation around $(X_1,Y_1)=(0,0)$ is given by
\begin{equation}
\begin{pmatrix}
a_1(t)(c_1-c_2)(c_1-c_3)\ldots(c_1-c_n)(c_1-t) & *\\
0 & \frac{a_1(t)(c_1-c_2)(c_1-c_3)\ldots(c_1-c_n)(c_1-t)}{m_1}
\end{pmatrix},
\end{equation}
where $* \in {\mathbb C}(t)$. We see that the local index at $(X_1,Y_1)=(0,0)$ is given by
\begin{equation}
\left(a_1(t)(c_1-c_2)(c_1-c_3)\ldots(c_1-c_n)(c_1-t),\frac{a_1(t)(c_1-c_2)(c_1-c_3)\ldots(c_1-c_n)(c_1-t)}{m_1} \right).
\end{equation}
The ratio of this local index is $m_1$.

For the remaining accessible singular points, we can discuss in the same way of this case.

In the coordinate system $(X_2,Y_2)=(z_1,1/w_1)$ we see that the system \eqref{fc2} admits $X_2=Y_2=0$ as its accessible singular points. The local index at $(X_2,Y_2)=(0,0)$ is given by
\begin{equation}
\left(-a_1(t),-n a_1(t)+\frac{a_1(t)}{m_1m_2 \ldots m_{n+1}}(m_1m_3 \ldots m_{n}m_{n+1}+m_2m_3 \ldots m_{n}m_{n+1}+\ldots m_1m_2 \ldots m_{n-1}m_{n}) \right).
\end{equation}
The ratio of this local index is given by
\begin{align}
\begin{split}
&\frac{-a_1(t)}{-n a_1(t)+\frac{a_1(t)}{m_1m_2 \ldots m_{n+1}}(m_1m_3 \ldots m_{n}m_{n+1}+m_2m_3 \ldots m_{n}m_{n+1}+\ldots +m_1m_2 \ldots m_{n-1}m_{n})}\\
&=\frac{-1}{-n+\frac{1}{m_1}+\frac{1}{m_2}+\ldots+\frac{1}{m_{n+1}}}.
\end{split}
\end{align}
Setting
\begin{equation}
m_{n+2}=\frac{-1}{-n+\frac{1}{m_1}+\frac{1}{m_2}+\ldots+\frac{1}{m_{n+1}}},
\end{equation}
we can obtain the relation \eqref{relation}.

Thus, we have completed the proof of Proposition \ref{prop1} and Theorem \ref{nonlinear:th}.

\newpage

\section{Appendix A}
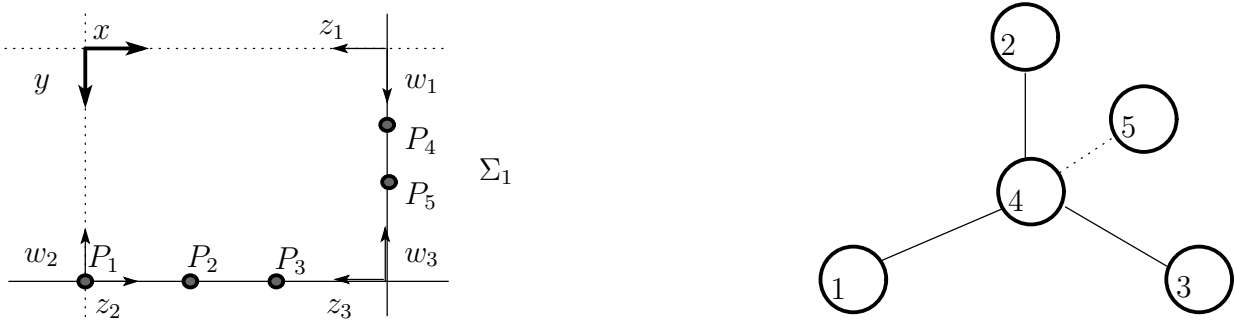
\begin{figure}[h]
\unitlength 0.1in
\begin{picture}( 64.2100, 16.4100)(  6.4000,-25.3000)
%
\special{pn 8}%
\special{pa 640 1120}%
\special{pa 2970 1120}%
\special{dt 0.045}%
%
\special{pn 8}%
\special{pa 660 2340}%
\special{pa 2960 2340}%
\special{fp}%
%
\special{pn 8}%
\special{pa 1060 940}%
\special{pa 1060 2530}%
\special{dt 0.045}%
%
\special{pn 8}%
\special{pa 2640 940}%
\special{pa 2640 2520}%
\special{fp}%
\put(31.2000,-18.3000){\makebox(0,0)[lb]{${\Sigma}_1$}}%
%
\special{pn 20}%
\special{pa 1060 1120}%
\special{pa 1360 1120}%
\special{fp}%
\special{sh 1}%
\special{pa 1360 1120}%
\special{pa 1294 1100}%
\special{pa 1308 1120}%
\special{pa 1294 1140}%
\special{pa 1360 1120}%
\special{fp}%
%
\special{pn 20}%
\special{pa 1060 1120}%
\special{pa 1060 1400}%
\special{fp}%
\special{sh 1}%
\special{pa 1060 1400}%
\special{pa 1080 1334}%
\special{pa 1060 1348}%
\special{pa 1040 1334}%
\special{pa 1060 1400}%
\special{fp}%
\put(11.0000,-10.7000){\makebox(0,0)[lb]{$x$}}%
\put(7.9000,-13.6000){\makebox(0,0)[lb]{$y$}}%
%
\special{pn 8}%
\special{pa 1060 2340}%
\special{pa 1060 2070}%
\special{fp}%
\special{sh 1}%
\special{pa 1060 2070}%
\special{pa 1040 2138}%
\special{pa 1060 2124}%
\special{pa 1080 2138}%
\special{pa 1060 2070}%
\special{fp}%
%
\special{pn 8}%
\special{pa 1050 2340}%
\special{pa 1330 2340}%
\special{fp}%
\special{sh 1}%
\special{pa 1330 2340}%
\special{pa 1264 2320}%
\special{pa 1278 2340}%
\special{pa 1264 2360}%
\special{pa 1330 2340}%
\special{fp}%
%
\special{pn 8}%
\special{pa 2640 1120}%
\special{pa 2640 1400}%
\special{fp}%
\special{sh 1}%
\special{pa 2640 1400}%
\special{pa 2660 1334}%
\special{pa 2640 1348}%
\special{pa 2620 1334}%
\special{pa 2640 1400}%
\special{fp}%
%
\special{pn 8}%
\special{pa 2630 1120}%
\special{pa 2360 1120}%
\special{fp}%
\special{sh 1}%
\special{pa 2360 1120}%
\special{pa 2428 1140}%
\special{pa 2414 1120}%
\special{pa 2428 1100}%
\special{pa 2360 1120}%
\special{fp}%
%
\special{pn 8}%
\special{pa 2630 2340}%
\special{pa 2630 2060}%
\special{fp}%
\special{sh 1}%
\special{pa 2630 2060}%
\special{pa 2610 2128}%
\special{pa 2630 2114}%
\special{pa 2650 2128}%
\special{pa 2630 2060}%
\special{fp}%
%
\special{pn 8}%
\special{pa 2620 2330}%
\special{pa 2370 2330}%
\special{fp}%
\special{sh 1}%
\special{pa 2370 2330}%
\special{pa 2438 2350}%
\special{pa 2424 2330}%
\special{pa 2438 2310}%
\special{pa 2370 2330}%
\special{fp}%
\put(22.8000,-10.7000){\makebox(0,0)[lb]{$z_1$}}%
\put(27.3000,-13.5000){\makebox(0,0)[lb]{$w_1$}}%
\put(7.4000,-22.6000){\makebox(0,0)[lb]{$w_2$}}%
\put(11.1000,-25.3000){\makebox(0,0)[lb]{$z_2$}}%
\put(27.3000,-22.6000){\makebox(0,0)[lb]{$w_3$}}%
\put(23.2000,-25.3000){\makebox(0,0)[lb]{$z_3$}}%
%
\special{pn 20}%
\special{sh 0.600}%
\special{ar 1060 2340 36 32  0.0000000 6.2831853}%
%
\special{pn 20}%
\special{sh 0.600}%
\special{ar 1610 2340 36 32  0.0000000 6.2831853}%
%
\special{pn 20}%
\special{sh 0.600}%
\special{ar 2060 2340 36 32  0.0000000 6.2831853}%
%
\special{pn 20}%
\special{sh 0.600}%
\special{ar 2640 1520 36 32  0.0000000 6.2831853}%
%
\special{pn 20}%
\special{sh 0.600}%
\special{ar 2650 1820 36 32  0.0000000 6.2831853}%
\put(15.9000,-22.9000){\makebox(0,0)[lb]{$P_2$}}%
\put(20.5000,-23.0000){\makebox(0,0)[lb]{$P_3$}}%
\put(10.7000,-22.9000){\makebox(0,0)[lb]{$P_1$}}%
\put(27.3000,-16.7000){\makebox(0,0)[lb]{$P_4$}}%
\put(27.3000,-19.6000){\makebox(0,0)[lb]{$P_5$}}%
%
\special{pn 20}%
\special{ar 5980 1060 172 172  0.0000000 6.2831853}%
%
\special{pn 20}%
\special{ar 5080 2330 172 172  0.0000000 6.2831853}%
%
\special{pn 20}%
\special{ar 6890 2330 172 172  0.0000000 6.2831853}%
%
\special{pn 20}%
\special{ar 6010 1870 172 172  0.0000000 6.2831853}%
%
\special{pn 20}%
\special{ar 6600 1490 172 172  0.0000000 6.2831853}%
%
\special{pn 8}%
\special{pa 5980 1250}%
\special{pa 5980 1690}%
\special{fp}%
%
\special{pn 8}%
\special{pa 5230 2230}%
\special{pa 5860 1960}%
\special{fp}%
%
\special{pn 8}%
\special{pa 6190 1950}%
\special{pa 6720 2260}%
\special{fp}%
\put(49.6000,-24.3000){\makebox(0,0)[lb]{$1$}}%
\put(58.5000,-11.5000){\makebox(0,0)[lb]{$2$}}%
\put(67.7000,-24.2000){\makebox(0,0)[lb]{$3$}}%
\put(58.9000,-19.7000){\makebox(0,0)[lb]{$4$}}%
\put(64.8000,-15.8000){\makebox(0,0)[lb]{$5$}}%
%
\special{pn 8}%
\special{pa 6160 1770}%
\special{pa 6450 1580}%
\special{dt 0.045}%
\end{picture}%
\label{fig:F1surfacePVI}
\caption{Hirzebruch surface ${\Sigma_1}$ (or ${\mathbb F}_1$ surface) }
\end{figure}

\begin{center}
\begin{tabular}{|c||c|} \hline
Equation & Painlev\'e VI system \eqref{PVI} with canonical Hamiltonian \eqref{HVI}  \\ \hline
Compactification & ${\Sigma}_1$ (or ${\mathbb F}_1$ surface)  \\ \hline
Accessible singular points & $(z_2,w_2)=(0,0),(1,0),(t,0), \ (z_1,w_1)=(0,-\alpha_2),(0,-(\alpha_2+\alpha_1))$  \\ \hline
Painlev\'e scheme & $
\begin{pmatrix}
P_1:(z_2,w_2)=(0,0), &P_2:(z_2,w_2)=(1,0), & P_3:(z_2,w_2)=(t,0)\\
\begin{pmatrix}
n_1 & \alpha_4\\
0 & 1 
\end{pmatrix} & \begin{pmatrix}
n_2 & \alpha_3\\
0 & 1 
\end{pmatrix} & \begin{pmatrix}
n_3 & \alpha_0\\
0 & 1 
\end{pmatrix}
\end{pmatrix}
$  \\ \hline
 & $
\begin{pmatrix}
P_4:(z_1,w_1)=(0,-\alpha_2), &P_5:(z_1,w_1)=(0,-(\alpha_2+\alpha_1))\\
\begin{pmatrix}
1 & 0\\
0 & n_4
\end{pmatrix} & \begin{pmatrix}
1 & 0\\
0 & n_5
\end{pmatrix}
\end{pmatrix}
$  \\ \hline
Relation of eigenvalues $n_i$ & $n_1 n_2 n_3(n_4 n_5+n_4+n_5)-(n_1 n_2+n_1 n_3+n_2 n_3)(n_4+n_5)=0$  \\ \hline
Painlev\'e VI case & $(n_1,n_2,n_3,n_4,n_5)=(2,2,2,1,1)$  \\ \hline
\end{tabular}
\end{center}
Here, we review the Hirzebruch surface ${\Sigma_1}$ (or ${\mathbb F}_1$ surface), which is obtained by gluing four copies of ${\mathbb C}^2$ via the following identification:
\begin{align}
\begin{split}
&U_j \cong {\mathbb C}^2 \ni (z_j,w_j) \ (j=0,1,2,3)\\
&z_0=x, \ w_0=y, \quad z_1=\frac{1}{x}, \ w_1=x y, \quad z_2=z_0, \ w_2=\frac{1}{w_0}, \quad z_3=z_1, \ w_3=\frac{1}{w_1}.
\end{split}
\end{align}
It is known that Painlev\'e VI system \eqref{PVI} with \eqref{HVI} in the coordinate system $(z_1,w_1)$ has two accessible singular points $P_4,P_5$. This differential system at each of accessible singular points $P_4,P_5$ passes the Painlev\'e $\alpha$-test (see \eqref{poiuy}).

We remark that we can not consider the case of compactification ${\mathbb P}^2$ because this system in the coordinate system $(X,Y)=(1/x,y/x)$ does not pass the Painlev\'e $\alpha$-method.

By a direct calculation, the above eigenvalue's relation can be transformed into the one obtained in Appendix B;
\begin{equation*}
\frac{1}{n_1}+\frac{1}{n_2}+\frac{1}{n_3}+\frac{1}{n_4}+\frac{1}{n_5}=3.
\end{equation*}

\section{Appendix B}

\begin{align*}
\begin{split}
&\frac{dx}{dt} =\frac{\partial \tilde{H}_{VI}}{\partial y}, \quad \frac{dy}{dt} =-\frac{\partial \tilde{H}_{VI}}{\partial x} \quad (\eta \in {\mathbb C}-\{0,1 \}),\\
&\tilde{H}_{VI}=\frac{1}{t(t-1)(t-\eta)}[x(x-1)(x-\eta)(x-t) y^2-\{ \alpha_1 (t-\eta)(x-1)x+2 \alpha_2 x(x-1)(x-\eta)\\
&+\alpha_3(t-1)(x-\eta)x+\alpha_4 t (x-1)(x-\eta)\}y+\alpha_2 \{(\alpha_1+\alpha_2)(t-\eta)+\alpha_2(x-1)+\alpha_3(t-1)+t \alpha_4 \} x].
\end{split}
\end{align*}

\begin{center}
\begin{tabular}{|c||c|} \hline
Equation & Painlev\'e VI system with the above symmetric Hamiltonian  \\ \hline
Compactification & ${\Sigma}_1$ (or ${\mathbb F}_1$ surface)  \\ \hline
Accessible singular points & $(z_2,w_2)=(0,0),(1,0),(t,0),(\eta,0), \ (z_1,w_1)=(0,-\alpha_2)$  \\ \hline
Painlev\'e scheme & $
\begin{pmatrix}
P_1:(z_2,w_2)=(0,0), &P_2:(z_2,w_2)=(1,0), & P_3:(z_2,w_2)=(t,0)\\
\begin{pmatrix}
n_4 & \alpha_4\\
0 & 1 
\end{pmatrix} & \begin{pmatrix}
n_3 & \alpha_3\\
0 & 1 
\end{pmatrix} & \begin{pmatrix}
n_0 & \alpha_0\\
0 & 1 
\end{pmatrix}
\end{pmatrix}
$  \\ \hline
 & $
\begin{pmatrix}
P_4:(z_2,w_2)=(\eta,0), &P_5:(z_1,w_1)=(0,-\alpha_2)\\
\begin{pmatrix}
n_1 & \alpha_1 \\
0 & 1
\end{pmatrix} & \begin{pmatrix}
1 & 0\\
* & n_2
\end{pmatrix}
\end{pmatrix}
$  \\ \hline
Relation of eigenvalues $n_i$ & $n_0 n_1 (n_2+1) n_3 n_4-n_0 n_1 n_3-n_0 n_1 n_4-n_0 n_3 n_4-n_1 n_3 n_4=0$  \\ \hline
Painlev\'e VI case & $(n_1,n_2,n_3,n_4,n_5)=(2,2,2,2,1)$  \\ \hline
\end{tabular}
\end{center}
It is known that symmetric Painlev\'e VI system at $P_5$ passes the Painlev\'e $\alpha$-method;
\begin{align}
\frac{d}{dT}\begin{pmatrix}
             Z \\
             W 
             \end{pmatrix}&=\begin{pmatrix}
             a-\frac{b W}{Z} \\
             \frac{a W}{Z}-\frac{b W^2}{Z^2}
             \end{pmatrix}=\frac{1}{Z^2} \begin{pmatrix}
             -b & 0  \\
             a & -b
             \end{pmatrix}\begin{pmatrix}
             Z W \\
             W^2 
             \end{pmatrix}+\begin{pmatrix}
             a \\
             0 
             \end{pmatrix} \ (a,b \in {\mathbb C}).
             \end{align}
Let us solve this system explicitly;
\begin{equation}
Z(T)=(a-b C_1)T+C_2, \quad W(T)=C_1 \{(a-b C_1)T+C_2 \} \quad (C_1,C_2 \in {\mathbb C}).
\end{equation}

In the case of $a=0$, let us consider a generalization of the above system;
\begin{align}
\frac{d}{dT}\begin{pmatrix}
             Z \\
             W 
             \end{pmatrix}&=\begin{pmatrix}
             -\frac{b W}{Z} \\
             -\frac{n_2 b W^2}{Z^2}
             \end{pmatrix}=-\frac{b}{Z^2} \begin{pmatrix}
             1 & 0  \\
             0 & n_2
             \end{pmatrix}\begin{pmatrix}
             Z W \\
             W^2 
             \end{pmatrix} \ (b \in {\mathbb C}).
             \end{align}
We solve this system explicitly;
\begin{equation*}
Z(T)=\{(n_2-2)(b C_1 T-C_2) \}^{\frac{1}{2-n_2}}, \quad W(T)=C_1 \{(n_2-2)(b C_1 T-C_2) \}^{\frac{n_2}{2-n_2}} \quad (C_1,C_2 \in {\mathbb C}).
\end{equation*}
Setting $\frac{1}{2-n_2}=N_2 \ (n_2=2-\frac{1}{N_2}, \ N_2 \in {\mathbb Z})$, we can obtain the following relation of eigenvalues $n_i \ (i=0,1,3,4)$ and $N_2$;
\begin{equation}
\frac{1}{n_0}+\frac{1}{n_1}+\frac{1}{N_2}+\frac{1}{n_3}+\frac{1}{n_4}=3.
\end{equation}
This equation has symmetry of symmetric group of degree five.

We remark that this relation coincides with the one obtained in the case of five accessible singular points in ${\Sigma}_3$.

\section{Appendix C}

{\bf Polynomial Hamiltonian of the sixth Painlev\'e system (see \cite{T1})}
\begin{align}\label{PVIC}
\begin{split}
&\frac{dq}{dt}=\frac{\partial H_{VI}}{\partial p}, \quad \frac{dp}{dt}=-\frac{\partial H_{VI}}{\partial q},\\
&H_{VI}(q,p,t;\alpha_0,\alpha_1,\alpha_2,\alpha_3,\alpha_4)\\
&=\frac{1}{t(t-1)}[p^2(q-t)(q-1)q-\{(\alpha_0-1)(q-1)q+\alpha_3(q-t)q\\
&+\alpha_4(q-t)(q-1)\}p+\alpha_2(\alpha_1+\alpha_2)(q-t)] \quad (\alpha_0+\alpha_1+2\alpha_2+\alpha_3+\alpha_4=1).
\end{split}
\end{align}

{\bf Holomorphy conditions (see \cite{T1})}
\begin{align}\label{holoPVI}
\begin{split}
&r_0:x_0=-((q-t)p-\alpha_0)p,\ y_0=\frac{1}{p}, \quad r_1:x_1=\frac{1}{q}, \ y_1=-(pq+\alpha_1+\alpha_2)q, \\
&r_2:x_2=\frac{1}{q}, \ y_2=-(pq+\alpha_2)q, \quad r_3:x_3=-((q-1)p-\alpha_3)p,\ y_3=\frac{1}{p},\\
&r_4:x_4=-(qp-\alpha_4)p,\ y_4=\frac{1}{p}.
\end{split}
\end{align}
Each transformation $r_i$ is birational and symplectic: $dy_i \wedge dx_i=dp \wedge dq$.

{\bf Symmetry}

The system \eqref{PVIC} is invariant under the following birational and symplectic transformations, whose generators $s_i \ (i=0,1,2,3,4)$ and ${\sigma}_i \ (i=1,2,3)$ are given by  (see \cite{Oka5})
\begin{align}\label{D4}
\begin{split}
s_0(q,p,t;\alpha_0,\alpha_1,\dots,\alpha_4) \rightarrow &(q,p-\frac{\alpha_0}{q-t},t;-\alpha_0,\alpha_1,\alpha_2+\alpha_0,\alpha_3,\alpha_4),\\
s_1(q,p,t;\alpha_0,\alpha_1,\dots,\alpha_4) \rightarrow &(q,p,t;\alpha_0,-\alpha_1,\alpha_2+\alpha_1,\alpha_3,\alpha_4),\\
s_2(q,p,t;\alpha_0,\alpha_1,\dots,\alpha_4) \rightarrow &(q+\frac{\alpha_2}{p},p,t;\alpha_0+\alpha_2,\alpha_1+\alpha_2,-\alpha_2,\\
&\alpha_3+\alpha_2,\alpha_4+\alpha_2),\\
s_3(q,p,t;\alpha_0,\alpha_1,\dots,\alpha_4) \rightarrow &(q,p-\frac{\alpha_3}{q-1},t;\alpha_0,\alpha_1,\alpha_2+\alpha_3,-\alpha_3,\alpha_4),\\
s_4(q,p,t;\alpha_0,\alpha_1,\dots,\alpha_4) \rightarrow &(q,p-\frac{\alpha_4}{q},t;\alpha_0,\alpha_1,\alpha_2+\alpha_4,\alpha_3,-\alpha_4),\\
{\sigma}_1(q,p,t;\alpha_0,\alpha_1,\dots,\alpha_4) \rightarrow &(1-q,-p,1-t;\alpha_0,\alpha_1,\alpha_2,\alpha_4,\alpha_3),\\
{\sigma}_2(q,p,t;\alpha_0,\alpha_1,\dots,\alpha_4) \rightarrow &(\frac{1}{q},-(pq+\alpha_2)q,\frac{1}{t};\alpha_0,\alpha_4,\alpha_2,\alpha_3,\alpha_1),\\
{\sigma}_3(q,p,t;\alpha_0,\alpha_1,\dots,\alpha_4) \rightarrow &(\frac{t-q}{t-1},-(t-1)p,\frac{t}{t-1};\alpha_4,\alpha_1,\alpha_2,\alpha_3,\alpha_0).
\end{split}
\end{align}
The system \eqref{PVIC} admits affine Weyl group symmetry of type $D_4^{(1)}$ as the group $<s_0,s_1,\ldots,s_4>$ of its B{\"a}cklund transformations.

Here, we review a relation between holomorphy and symmetry conditions of Painlev\'e type systems. For example, we will consider the following relations (see Figure 5):
$$
(X,Y)=\left(\frac{1}{q},-(p q+\alpha_2)q \right) \Longleftrightarrow (X',Y')=\left(q+\frac{\alpha_2}{p},p \right),
$$
and
$$
(X,Y)=\left(-(q p-\alpha_i)p,\frac{1}{p} \right) \Longleftrightarrow (X',Y')=\left(q,p-\frac{\alpha_i}{q} \right).
$$

\begin{figure}
\unitlength 0.1in
\begin{picture}( 52.7000, 36.6000)(  8.5000,-41.0000)
%
\special{pn 8}%
\special{pa 1080 940}%
\special{pa 2770 940}%
\special{fp}%
\special{sh 1}%
\special{pa 2770 940}%
\special{pa 2704 920}%
\special{pa 2718 940}%
\special{pa 2704 960}%
\special{pa 2770 940}%
\special{fp}%
%
\special{pn 8}%
\special{pa 1790 590}%
\special{pa 1220 1500}%
\special{fp}%
\special{sh 1}%
\special{pa 1220 1500}%
\special{pa 1272 1454}%
\special{pa 1248 1456}%
\special{pa 1238 1434}%
\special{pa 1220 1500}%
\special{fp}%
\put(8.5000,-16.8000){\makebox(0,0)[lb]{$\frac{1}{x}$}}%
\put(22.2000,-8.8000){\makebox(0,0)[lb]{$\frac{y}{x}$}}%
%
\special{pn 20}%
\special{pa 1990 2190}%
\special{pa 1990 1830}%
\special{fp}%
\special{sh 1}%
\special{pa 1990 1830}%
\special{pa 1970 1898}%
\special{pa 1990 1884}%
\special{pa 2010 1898}%
\special{pa 1990 1830}%
\special{fp}%
%
\special{pn 8}%
\special{pa 1060 2660}%
\special{pa 2750 2660}%
\special{fp}%
\special{sh 1}%
\special{pa 2750 2660}%
\special{pa 2684 2640}%
\special{pa 2698 2660}%
\special{pa 2684 2680}%
\special{pa 2750 2660}%
\special{fp}%
%
\special{pn 8}%
\special{pa 1780 2390}%
\special{pa 1210 3540}%
\special{dt 0.045}%
%
\special{pn 8}%
\special{pa 1210 3220}%
\special{pa 1660 4100}%
\special{fp}%
\special{sh 1}%
\special{pa 1660 4100}%
\special{pa 1648 4032}%
\special{pa 1636 4054}%
\special{pa 1612 4050}%
\special{pa 1660 4100}%
\special{fp}%
\put(22.1000,-26.0000){\makebox(0,0)[lb]{$\frac{y}{x}$}}%
\put(14.6000,-41.1000){\makebox(0,0)[lb]{$\frac{1}{x}$}}%
%
\special{pn 20}%
\special{pa 1290 3380}%
\special{pa 1430 3090}%
\special{fp}%
\special{sh 1}%
\special{pa 1430 3090}%
\special{pa 1384 3142}%
\special{pa 1408 3138}%
\special{pa 1420 3160}%
\special{pa 1430 3090}%
\special{fp}%
\put(15.0000,-31.9000){\makebox(0,0)[lb]{$y$}}%
%
\special{pn 20}%
\special{pa 4320 800}%
\special{pa 4320 440}%
\special{fp}%
\special{sh 1}%
\special{pa 4320 440}%
\special{pa 4300 508}%
\special{pa 4320 494}%
\special{pa 4340 508}%
\special{pa 4320 440}%
\special{fp}%
%
\special{pn 8}%
\special{pa 3526 1188}%
\special{pa 6120 1188}%
\special{fp}%
\special{sh 1}%
\special{pa 6120 1188}%
\special{pa 6054 1168}%
\special{pa 6068 1188}%
\special{pa 6054 1208}%
\special{pa 6120 1188}%
\special{fp}%
%
\special{pn 8}%
\special{pa 4340 820}%
\special{pa 3572 2250}%
\special{fp}%
%
\special{pn 8}%
\special{pa 3510 1676}%
\special{pa 4110 3474}%
\special{dt 0.045}%
%
\special{pn 8}%
\special{pa 3602 3104}%
\special{pa 5830 3856}%
\special{fp}%
\special{sh 1}%
\special{pa 5830 3856}%
\special{pa 5772 3816}%
\special{pa 5778 3840}%
\special{pa 5760 3854}%
\special{pa 5830 3856}%
\special{fp}%
\put(53.2000,-39.6000){\makebox(0,0)[lb]{$\frac{1}{x}$}}%
%
\special{pn 20}%
\special{pa 3648 2102}%
\special{pa 3940 1572}%
\special{fp}%
\special{sh 1}%
\special{pa 3940 1572}%
\special{pa 3890 1622}%
\special{pa 3914 1620}%
\special{pa 3926 1640}%
\special{pa 3940 1572}%
\special{fp}%
%
\special{pn 20}%
\special{pa 3634 2088}%
\special{pa 3802 2516}%
\special{fp}%
\special{sh 1}%
\special{pa 3802 2516}%
\special{pa 3796 2446}%
\special{pa 3782 2466}%
\special{pa 3760 2460}%
\special{pa 3802 2516}%
\special{fp}%
%
\special{pn 20}%
\special{pa 4032 3252}%
\special{pa 3880 2780}%
\special{fp}%
\special{sh 1}%
\special{pa 3880 2780}%
\special{pa 3882 2850}%
\special{pa 3896 2832}%
\special{pa 3920 2838}%
\special{pa 3880 2780}%
\special{fp}%
\put(40.1700,-17.3400){\makebox(0,0)[lb]{$y$}}%
\put(38.1700,-24.7100){\makebox(0,0)[lb]{$\frac{1}{xy}$}}%
\put(39.5500,-28.3900){\makebox(0,0)[lb]{$xy$}}%
%
\special{pn 20}%
\special{sh 0.600}%
\special{ar 3830 2640 16 30  0.0000000 6.2831853}%
%
\special{pn 20}%
\special{sh 0.600}%
\special{ar 1290 3380 16 30  0.0000000 6.2831853}%
\put(31.0000,-21.8000){\makebox(0,0)[lb]{$U_2$}}%
\put(35.8000,-34.4000){\makebox(0,0)[lb]{$U_1$}}%
%
\special{pn 20}%
\special{sh 0.600}%
\special{ar 1580 950 16 30  0.0000000 6.2831853}%
\put(14.1000,-21.4000){\makebox(0,0)[lb]{Step 1}}%
\put(37.5000,-7.5000){\makebox(0,0)[lb]{Step 2}}%
\put(35.2000,-27.2000){\makebox(0,0)[lb]{$P$}}%
%
\special{pn 13}%
\special{pa 4410 2620}%
\special{pa 3940 2620}%
\special{fp}%
\special{sh 1}%
\special{pa 3940 2620}%
\special{pa 4008 2640}%
\special{pa 3994 2620}%
\special{pa 4008 2600}%
\special{pa 3940 2620}%
\special{fp}%
\put(44.3000,-26.8000){\makebox(0,0)[lb]{accessible singular point}}%
\put(20.8000,-21.5000){\makebox(0,0)[lb]{Blow up}}%
\put(44.5000,-7.5000){\makebox(0,0)[lb]{Blow up}}%
\end{picture}%
\label{Holomorphy}
\caption{After Step 2, let us take the coordinate neighborhood $\{U_1,(x_1,y_1)=(\frac{1}{q},p q)\}$. We see that the point $P:\{(x_1,y_1)=(0,-\alpha_2) \}$ is its accessible singular point, where the parameter $\alpha_2 \in {\mathbb C}$. Blowing up this point $P$, we can obtain the holomorphy condition $(X,Y)=(\frac{1}{q},-(p q+\alpha_2)q)$ (see \cite{T1}). In another coordinate neighborhood $\{U_2,(x_2,y_2)=(\frac{1}{q p},p)\}$, we easily see that the point $P:\{(x_2,y_2)=(-\frac{1}{\alpha_2},0)\}$ is an accessible singular point. Since $\alpha_2$ is not $0$, then we can replace the coordinate system $(x_2,y_2)$ as $(x'_2,y'_2)=(q p,p)$. In new coordinate system, we see that the point $P$ is given by $\{(x'_2,y'_2)=(-\alpha_2,0) \}$. Blowing up this point $P$, we can obtain the symmetry condition $(X',Y')=(q+\frac{\alpha_2}{p},p)$. Thus, both relations can be obtained by blowing up the accessible singular point $P$, respectively. We remark that in \cite{T1}  the holomorphy condition $(X,Y)=\left(-(q p-\alpha_i)p,\frac{1}{p} \right)$ was explained. This case is explained by similar way.
}
\end{figure}
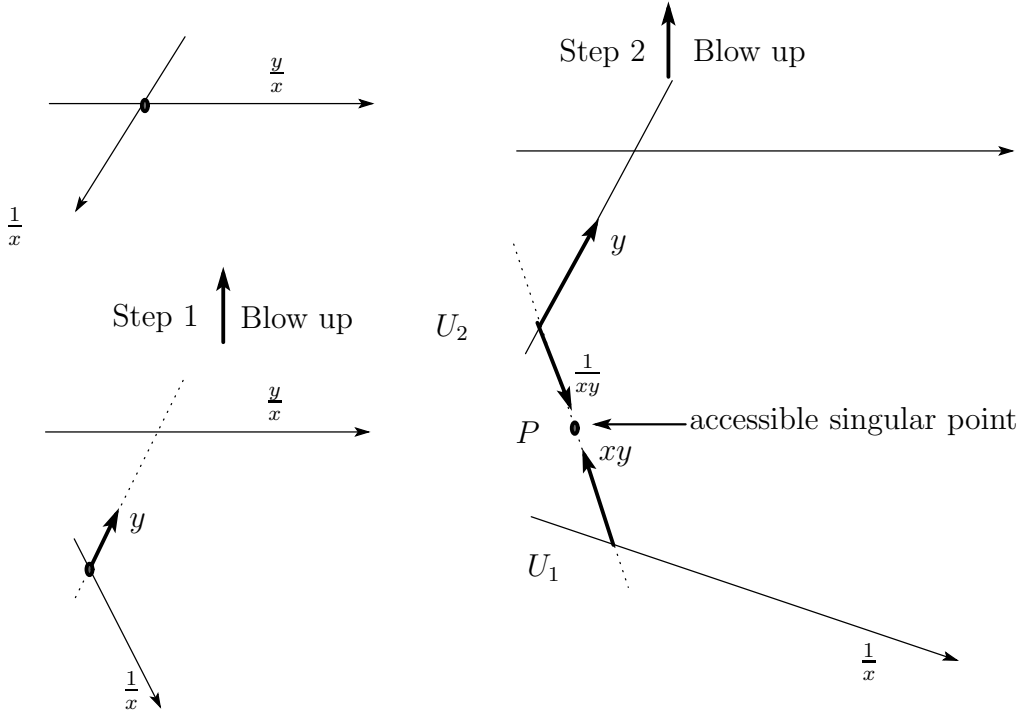

In this appendix, we study the polynomial Hamiltonian system (see \cite{T1,Yamada2}) given by
\begin{align}\label{SPVI}
\begin{split}
&\frac{dq}{dt}=\frac{\partial \tilde{H}_{VI}}{\partial p}, \quad \frac{dp}{dt}=-\frac{\partial \tilde{H}_{VI}}{\partial q},\\
&\tilde{H}_{VI}(q,p,t;\alpha_0,\alpha_1,\alpha_2,\alpha_3,\alpha_4)\\
&=-\frac{q^3p^4}{t(t-1)}-\frac{(\alpha_0+\alpha_3-2\alpha_4-1)q^2p^3}{t(t-1)}\\
&-\frac{(t+1)q^2p^2}{t(t-1)}-\frac{(\alpha_1\alpha_2+\alpha_2^2+2\alpha_4-2\alpha_0\alpha_4-2\alpha_3\alpha_4+\alpha_4^2)qp^2}{t(t-1)}\\
&-\frac{\{(\alpha_3-\alpha_4)t+\alpha_0-\alpha_4-1\}qp}{t(t-1)}-\frac{q}{t-1}+\frac{\alpha_4(\alpha_2+\alpha_4)(\alpha_1+\alpha_2+\alpha_4)p}{t(t-1)},
\end{split}
\end{align}
where $\tilde{H}_{VI}=r_4^{-1} (H_{VI})$. Here, we remark that the system \eqref{PVIC} is not invariant under the birational transformation $r_4$. This transformation is called {\it holomorphy} (see \cite{Yamada}).

We remark that we will see that this system has a 1-parameter family of formal Laurent series:
\begin{align}\label{Laurent series11}
\begin{split}
&q=-\frac{t_0}{T^2}+\frac{(1+\alpha_0-\alpha_3)t_0}{(t_0-1)T}+h+{\mathcal O}(T),\\
&p=T+\frac{\{(2+2\alpha_0-\alpha_3-\alpha_4)t_0+\alpha_4-\alpha_0-1 \}}{2t_0(t_0-1)}T^2+{\mathcal O}(T^3)
\end{split}
\end{align}
and
\begin{align}\label{Laurent series12}
\begin{split}
&q=-\frac{t_0^2}{T^2}-\frac{(t_0+\alpha_0-\alpha_3-2)t_0}{(t_0-1)T}+h+{\mathcal O}(T),\\
&p=-\frac{1}{t_0}T+\frac{\{(1-\alpha_3-\alpha_4)t_0+\alpha_0+\alpha_4-2\alpha_3-2 \}}{2t_0^2(t_0-1)}T^2+{\mathcal O}(T^3),
\end{split}
\end{align}
where $T:=t-t_0$, $h$ is its free parameter and  the symbol ${\mathcal O}$ denotes Landau symbol.

We will show that these Laurent series are convergent by using the holomorphy $\tilde{r}_0,\tilde{r}_3$, respectively;

\begin{align}
\begin{split}
&\tilde{r}_0:x_0=q+\frac{\alpha_0-\alpha_4}{p}+\frac{t}{p^2},\ y_0=p,\\
&\tilde{r}_3:x_3=q+\frac{\alpha_3-\alpha_4}{p}+\frac{1}{p^2},\ y_3=p.
\end{split}
\end{align}

By the transformation $r_4$ these Laurent series are transformed into the following series (See \cite{Oka2}; P 212):
\begin{align}\label{Laurent series3}
\begin{split}
&r_4(q)=t_0+(\alpha_0+1)T+h T^2+{\mathcal O}(T^3),\\
&r_4(p)=\frac{1}{T}(1+{\mathcal O}(T))
\end{split}
\end{align}
and
\begin{align}\label{Laurent series4}
\begin{split}
&r_4(q)=1-\frac{\alpha_3}{t_0}T+h T^2+{\mathcal O}(T^3),\\
&r_4(p)=-\frac{t_0}{T}(1+{\mathcal O}(T)).
\end{split}
\end{align}

These meromorphic solutions \eqref{Laurent series11},\eqref{Laurent series12} can be characterized by diagram automorphisms ${\sigma}_1,{\sigma}_3$ (see \eqref{D4}), respectively.

{\bf Holomorphy conditions}

\begin{align}\label{holoPVI2}
\begin{split}
&\tilde{r}_0:x_0=q+\frac{\alpha_0-\alpha_4}{p}+\frac{t}{p^2},\ y_0=p,\\
&r_1:x_1=-(qp-(\alpha_1+\alpha_2+\alpha_4))p, \ y_1=\frac{1}{p}, \\
&r_2:x_2=-(qp-(\alpha_2+\alpha_4))p, \ y_2=\frac{1}{p},\\
&\tilde{r}_3:x_3=q+\frac{\alpha_3-\alpha_4}{p}+\frac{1}{p^2},\ y_3=p,\\
&r_4:x_4=-(qp-\alpha_4)p,\ y_4=\frac{1}{p}.
\end{split}
\end{align}
The transformations $\tilde{r}_0,\tilde{r}_3$ were known as one of patching data of the third Painlev\'e system. These transformations can be constructed by successive blowing-up procedures of double accessible singular point (see \cite{MMT}).

It is still an open question whether we can obtain the Hamiltonian system \eqref{SPVI} by solving $3 \times 3$ Mazzocco's Lax pair  (cf. \cite{Mazzocco}) satisfying the following Riemann scheme:
\begin{equation}
\begin{pmatrix}
X=0 \ (2) & X=\infty\\
$\(\overbrace{\begin{matrix}
0\\
t \\
1
\end{matrix} \quad \begin{matrix}
0\\
\alpha_0-\alpha_4\\
\alpha_3-\alpha_4
\end{matrix} }\)$  & \begin{matrix}
\alpha_4\\
\alpha_4+\alpha_2 \\
\alpha_4+\alpha_2+\alpha_1
\end{matrix}
\end{pmatrix}
\end{equation}

We remark that the system \eqref{SPVI} is invariant under the following birational and symplectic transformations, whose generators $s_i \ (i=1,2,4)$ and ${\sigma}_i \ (i=1,3)$, are given by

{\bf Symmetry}
\begin{align*}\label{F4}
\begin{split}
{\sigma}_3(q,p,t;\alpha_0,\alpha_1,\dots,\alpha_4) \rightarrow &((1-t)\left(q+\frac{\alpha_0-\alpha_4}{p}+\frac{t}{p^2}\right),\frac{p}{1-t},\frac{t}{t-1};\\
&\alpha_4,\alpha_1,\alpha_2,\alpha_3,\alpha_0),\\
s_1(q,p,t;\alpha_0,\alpha_1,\dots,\alpha_4) \rightarrow &(q,p,t;\alpha_0,-\alpha_1,\alpha_2+\alpha_1,\alpha_3,\alpha_4),\\
s_2(q,p,t;\alpha_0,\alpha_1,\dots,\alpha_4) \rightarrow &(q,p,t;\alpha_0+\alpha_2,\alpha_1+\alpha_2,-\alpha_2,\\
&\alpha_3+\alpha_2,\alpha_4+\alpha_2),\\
{\sigma}_1(q,p,t;\alpha_0,\alpha_1,\dots,\alpha_4) \rightarrow &(-\left(q+\frac{\alpha_3-\alpha_4}{p}+\frac{1}{p^2}\right),-p,1-t;\\
&\alpha_0,\alpha_1,\alpha_2,\alpha_4,\alpha_3),\\
s_4(q,p,t;\alpha_0,\alpha_1,\dots,\alpha_4) \rightarrow &\left(q,p-\frac{\alpha_4}{q},t;\alpha_0,\alpha_1,\alpha_2+\alpha_4,\alpha_3,-\alpha_4 \right).
\end{split}
\end{align*}

\begin{figure}
\unitlength 0.1in
\begin{picture}(44.75,44.76)(14.00,-45.86)
%
\special{pn 8}%
\special{pa 1677 323}%
\special{pa 3182 323}%
\special{fp}%
%
\special{pn 8}%
\special{pa 1875 178}%
\special{pa 2675 1128}%
\special{fp}%
%
\special{pn 8}%
\special{pa 2992 194}%
\special{pa 2216 1098}%
\special{fp}%
%
\special{pn 8}%
\special{pa 2454 1524}%
\special{pa 1670 2101}%
\special{fp}%
\special{pa 2207 1546}%
\special{pa 3048 2094}%
\special{fp}%
%
\special{pn 8}%
\special{pa 1724 1942}%
\special{pa 2445 2694}%
\special{fp}%
%
\special{pn 8}%
\special{pa 2999 1957}%
\special{pa 2200 2709}%
\special{fp}%
%
\special{pn 20}%
\special{pa 2683 3059}%
\special{pa 1709 3469}%
\special{fp}%
%
\special{pn 20}%
\special{pa 1749 3150}%
\special{pa 1962 4229}%
\special{fp}%
%
\special{pn 20}%
\special{pa 1756 4009}%
\special{pa 2643 4275}%
\special{fp}%
%
\special{pn 20}%
\special{pa 1717 3667}%
\special{pa 2493 3507}%
\special{fp}%
%
\special{pn 20}%
\special{pa 1724 3880}%
\special{pa 2660 3880}%
\special{fp}%
%
\special{pn 13}%
\special{pa 2375 2968}%
\special{pa 3427 3591}%
\special{fp}%
%
\special{pn 20}%
\special{pa 3237 3142}%
\special{pa 3523 4381}%
\special{fp}%
%
\special{pn 20}%
\special{pa 4663 2215}%
\special{pa 3673 3454}%
\special{fp}%
%
\special{pn 20}%
\special{pa 3752 2884}%
\special{pa 3966 3948}%
\special{fp}%
%
\special{pn 20}%
\special{pa 3966 2740}%
\special{pa 4251 3659}%
\special{fp}%
%
\special{pn 20}%
\special{pa 4227 2405}%
\special{pa 4663 3310}%
\special{fp}%
%
\special{pn 20}%
\special{pa 4425 2170}%
\special{pa 5875 3781}%
\special{fp}%
%
\special{pn 8}%
\special{pa 2224 3652}%
\special{pa 2493 3256}%
\special{dt 0.045}%
\special{pa 2493 3256}%
\special{pa 2493 3257}%
\special{dt 0.045}%
%
\special{pn 8}%
\special{pa 2326 4009}%
\special{pa 2739 3629}%
\special{dt 0.045}%
\special{pa 2739 3629}%
\special{pa 2739 3629}%
\special{dt 0.045}%
%
\special{pn 8}%
\special{pa 3514 3788}%
\special{pa 2984 4184}%
\special{dt 0.045}%
\special{pa 2984 4184}%
\special{pa 2985 4184}%
\special{dt 0.045}%
%
\special{pn 8}%
\special{pa 3625 4016}%
\special{pa 2960 4427}%
\special{dt 0.045}%
\special{pa 2960 4427}%
\special{pa 2961 4427}%
\special{dt 0.045}%
%
\special{pn 8}%
\special{pa 3610 4206}%
\special{pa 3142 4541}%
\special{dt 0.045}%
\special{pa 3142 4541}%
\special{pa 3143 4541}%
\special{dt 0.045}%
%
\special{pn 8}%
\special{pa 4259 3272}%
\special{pa 4069 3834}%
\special{dt 0.045}%
\special{pa 4069 3834}%
\special{pa 4069 3833}%
\special{dt 0.045}%
%
\special{pn 8}%
\special{pa 4647 2899}%
\special{pa 4497 3606}%
\special{dt 0.045}%
\special{pa 4497 3606}%
\special{pa 4497 3605}%
\special{dt 0.045}%
%
\special{pn 8}%
\special{pa 5234 2656}%
\special{pa 5043 3766}%
\special{dt 0.045}%
\special{pa 5043 3766}%
\special{pa 5043 3765}%
\special{dt 0.045}%
%
\special{pn 8}%
\special{pa 5573 2892}%
\special{pa 5352 4054}%
\special{dt 0.045}%
\special{pa 5352 4054}%
\special{pa 5352 4053}%
\special{dt 0.045}%
%
\special{pn 8}%
\special{pa 5828 3272}%
\special{pa 5677 4586}%
\special{dt 0.045}%
\special{pa 5677 4586}%
\special{pa 5677 4585}%
\special{dt 0.045}%
%
\special{pn 20}%
\special{sh 0.600}%
\special{ar 2881 323 17 24  0.0000000 6.2831853}%
%
\special{pn 20}%
\special{sh 0.600}%
\special{ar 1994 323 16 24  0.0000000 6.2831853}%
%
\special{pn 20}%
\special{sh 0.600}%
\special{ar 1796 2010 17 24  0.0000000 6.2831853}%
%
\special{pn 20}%
\special{sh 0.600}%
\special{ar 2937 2018 16 23  0.0000000 6.2831853}%
\put(47.9000,-25.7200){\makebox(0,0)[lb]{$D_4^{(1)}$-lattice}}%
\put(40.2900,-38.9500){\makebox(0,0)[lb]{$r_0$}}%
\put(44.3400,-36.1400){\makebox(0,0)[lb]{$r_3$}}%
\put(49.7900,-37.5000){\makebox(0,0)[lb]{$r_1$}}%
\put(52.8100,-40.7000){\makebox(0,0)[lb]{$r_2$}}%
\put(56.3000,-45.6400){\makebox(0,0)[lb]{$r_4$}}%
%
\special{pn 20}%
\special{pa 3863 612}%
\special{pa 5439 612}%
\special{fp}%
%
\special{pn 20}%
\special{pa 4244 505}%
\special{pa 3697 1501}%
\special{fp}%
%
\special{pn 20}%
\special{pa 5170 490}%
\special{pa 5828 1653}%
\special{fp}%
%
\special{pn 20}%
\special{pa 3642 1189}%
\special{pa 4038 1759}%
\special{fp}%
%
\special{pn 20}%
\special{pa 3744 1037}%
\special{pa 4315 1592}%
\special{fp}%
%
\special{pn 20}%
\special{pa 3887 870}%
\special{pa 4536 1265}%
\special{fp}%
%
\special{pn 8}%
\special{pa 4219 1379}%
\special{pa 4157 1881}%
\special{dt 0.045}%
\special{pa 4157 1881}%
\special{pa 4157 1880}%
\special{dt 0.045}%
%
\special{pn 8}%
\special{pa 4442 1030}%
\special{pa 4449 1546}%
\special{dt 0.045}%
\special{pa 4449 1546}%
\special{pa 4449 1545}%
\special{dt 0.045}%
%
\special{pn 8}%
\special{pa 5613 992}%
\special{pa 5091 1425}%
\special{dt 0.045}%
\special{pa 5091 1425}%
\special{pa 5092 1425}%
\special{dt 0.045}%
%
\special{pn 8}%
\special{pa 5764 1204}%
\special{pa 5177 1683}%
\special{dt 0.045}%
\special{pa 5177 1683}%
\special{pa 5178 1683}%
\special{dt 0.045}%
%
\special{pn 8}%
\special{pa 5820 1425}%
\special{pa 5328 1964}%
\special{dt 0.045}%
\special{pa 5328 1964}%
\special{pa 5328 1964}%
\special{dt 0.045}%
%
\special{pn 20}%
\special{pa 2437 2907}%
\special{pa 2437 2732}%
\special{fp}%
\special{sh 1}%
\special{pa 2437 2732}%
\special{pa 2417 2799}%
\special{pa 2437 2785}%
\special{pa 2457 2799}%
\special{pa 2437 2732}%
\special{fp}%
%
\special{pn 20}%
\special{pa 4600 110}%
\special{pa 4600 429}%
\special{fp}%
\special{sh 1}%
\special{pa 4600 429}%
\special{pa 4620 362}%
\special{pa 4600 376}%
\special{pa 4580 362}%
\special{pa 4600 429}%
\special{fp}%
%
\special{pn 20}%
\special{pa 4583 1812}%
\special{pa 4583 2162}%
\special{fp}%
\special{sh 1}%
\special{pa 4583 2162}%
\special{pa 4603 2095}%
\special{pa 4583 2109}%
\special{pa 4563 2095}%
\special{pa 4583 2162}%
\special{fp}%
%
\special{pn 8}%
\special{pa 2398 1121}%
\special{pa 2398 1432}%
\special{dt 0.045}%
\special{sh 1}%
\special{pa 2398 1432}%
\special{pa 2418 1365}%
\special{pa 2398 1379}%
\special{pa 2378 1365}%
\special{pa 2398 1432}%
\special{fp}%
\put(14.7200,-9.1600){\makebox(0,0)[lb]{${\mathbb P}^2$}}%
\put(14.0000,-25.8000){\makebox(0,0)[lb]{${\mathbb P}^1 \times {\mathbb P}^1$}}%
\put(20.8900,-32.1100){\makebox(0,0)[lb]{(-2)}}%
\put(31.3500,-31.3500){\makebox(0,0)[lb]{(-4)}}%
\put(38.5500,-7.9400){\makebox(0,0)[lb]{(-3)}}%
\put(53.3600,-8.0200){\makebox(0,0)[lb]{(-3)}}%
\put(27.9400,-32.1100){\makebox(0,0)[lb]{(-1)}}%
\put(45.7600,-5.8900){\makebox(0,0)[lb]{(-1)}}%
\put(15.9800,-31.2700){\makebox(0,0)[lb]{(-3)}}%
\end{picture}%
\label{fig:}
\caption{The space of initial conditions of \eqref{SPVI} can be constructed by the above way. The symbol $\bullet$ denotes its accessible singular points of the system \eqref{SPVI}. Each symbol $(*)$ denotes self-intersection number of each ${\mathbb P}^1$.}
\end{figure}
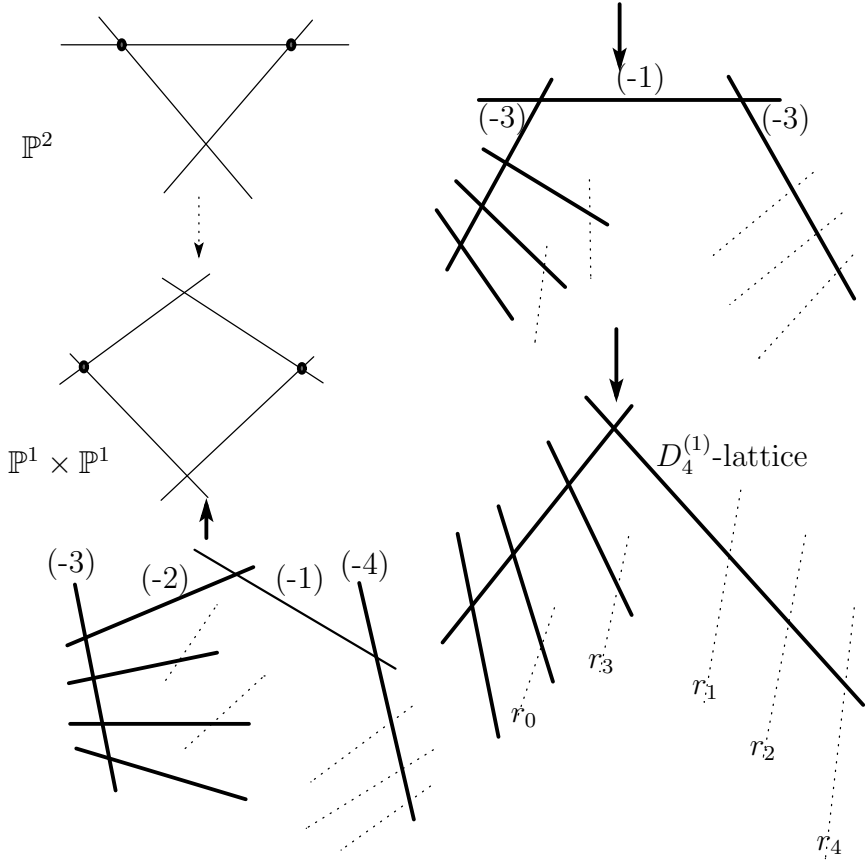

\begin{center}
\begin{tabular}{|c||c|} \hline
Equation & Painlev\'e VI system with another Hamiltonian \eqref{SPVI}  \\ \hline
Compact. & ${\Sigma}_1$ (or ${\mathbb F}_1$ surface)  \\ \hline
Accessible sing. & $(z_1,w_1)=(\alpha_4,0),(\alpha_2+\alpha_4,0),(\alpha_1+\alpha_2+\alpha_4,0),$  \\ \hline
 & $(Z_2,W_2)=(-1,0),(-t,0)$  \\ \hline
Painlev\'e scheme & $
\begin{pmatrix}
P_1:z_1=\alpha_4, &P_2:z_1=\alpha_2+\alpha_4, & P_3:z_1=\alpha_1+\alpha_2+\alpha_4\\
\begin{pmatrix}
n_1 & 0\\
0 & 1 
\end{pmatrix} & \begin{pmatrix}
n_2 & 0\\
0 & 1 
\end{pmatrix} & \begin{pmatrix}
n_3 & 0\\
0 & 1 
\end{pmatrix}
\end{pmatrix}
$  \\ \hline
 & $
\begin{pmatrix}
P_4:(Z_2,W_2)=(-1,0), &P_5:(Z_2,W_2)=(-t,0)\\
\begin{pmatrix}
n_4 & \alpha_4-\alpha_3\\
0 & 1
\end{pmatrix} & \begin{pmatrix}
n_5 & \alpha_4-\alpha_0\\
0 & 1
\end{pmatrix}
\end{pmatrix}
$, where $(Z_2,W_2)=(q p^2,p)$  \\ \hline
Relation of $n_i$ & $n_1 n_2 n_3 n_4 n_5+2n_1 n_2 n_3(n_4+n_5)+(n_1 n_2+n_1 n_3+n_2 n_3)(n_4 n_5-2n_4-2n_5)=0$  \\ \hline
Painlev\'e VI case & $(n_1,n_2,n_3,n_4,n_5)=(1,1,1,2,2)$  \\ \hline
\end{tabular}
\end{center}
Here, Hirzebruch surface ${\Sigma_1}$ (or ${\mathbb F}_1$ surface) is obtained by gluing four copies of ${\mathbb C}^2$ via the following identification:
\begin{align}
\begin{split}
&U_j \cong {\mathbb C}^2 \ni (z_j,w_j) \ (j=0,1,2,3)\\
&z_0=q, \ w_0=p, \quad z_1=q p, \ w_1=\frac{1}{p}, \quad z_2=\frac{1}{z_0}, \ w_2=w_0, \quad z_3=\frac{1}{z_1}, \ w_3=w_1.
\end{split}
\end{align}
Here, $(Z_2,W_2)=(q p^2,p)$. This coordinate system can be obtained by resolving a double accessible singular point (see Painlev\'e V case, \cite{MMT}).

By a direct calculation, the above eigenvalue's relation can be transformed into the one obtained in Appendix B;
\begin{equation*}
\frac{1}{n_1}+\frac{1}{n_2}+\frac{1}{n_3}+\frac{1}{n_4}+\frac{1}{n_5}=3.
\end{equation*}

\section{Appendix D}

The eigenvalue's relation;
\begin{equation}
\frac{1}{n_1}+\frac{1}{n_2}+\frac{1}{n_3}+\frac{1}{n_4}+\frac{1}{n_5}=b \quad (b \in {\mathbb C})
\end{equation}
can be transformed into a one-parameter family of quintic hypersurfaces (see \cite{Candelas});
\begin{equation}
x_1^5+x_2^5+x_3^5+x_4^5+x_5^5-b x_1 x_2 x_3 x_4 x_5=0,
\end{equation}
where we can make a change of variables:
\begin{align}
\begin{split}
&n_1=\frac{x_2 x_3 x_4 x_5}{x_1^4}, \quad n_2=\frac{x_1 x_3 x_4 x_5}{x_2^4}, \quad n_3=\frac{x_1 x_2 x_4 x_5}{x_3^4}, \quad n_4=-\frac{x_1 x_2 x_3 x_4 x_5 n_5}{x_1 x_2 x_3 x_4 x_5-(x_4^5+x_5^5) n_5}.
\end{split}
\end{align}

In general, the eigenvalue's relation (cf. \eqref{relation});
\begin{equation}
\frac{1}{n_1}+\frac{1}{n_2}+\cdots+\frac{1}{n_N}=b \quad (b \in {\mathbb C})
\end{equation}
can be transformed into a one-parameter family of Calabi-Yau hypersurfaces;
\begin{equation}
x_1^N+x_2^N+\cdots+x_N^N-b x_1 x_2 \cdots x_N=0,
\end{equation}
where we can make a change of variables:
\begin{align}
\begin{split}
&n_1=\frac{x_2 x_3 \cdots x_N}{x_1^{N-1}}, \quad n_2=\frac{x_1 x_3 \cdots x_N}{x_2^{N-1}}, \cdots ,n_{N-2}=\frac{x_1 \cdots x_{N-3} x_{N-1} x_N}{x_{N-2}^{N-1}},\\
&n_{N-1}=-\frac{x_1 x_2 \cdots x_N n_N}{x_1 x_2 \cdots x_N-(x_{N-1}^N+x_N^N) n_N}.
\end{split}
\end{align}

We remark that the eigenvalue's relation;
\begin{equation}\label{n3}
\frac{1}{n_1}+\frac{1}{n_2}+\frac{1}{n_3}=b \quad (b \in {\mathbb C})
\end{equation}
can be transformed into a cubic surface (cf. \cite{Iwasaki1,Dubrovin1});
\begin{equation}
b(x_1 x_2 x_3+x_1^2+x_2^2+x_3^2)=0,
\end{equation}
where we can make a change of variables:
\begin{align}
\begin{split}
&n_1=-\frac{x_2 x_3}{b x_1}, \quad n_2=-\frac{x_1 x_3}{b x_2}, \quad n_{3}=-\frac{x_1 x_2}{b x_3}.
\end{split}
\end{align}
On the other hand, the eigenvalue's relation \eqref{n3} can be transformed into a one-parameter family of cubic surfaces;
\begin{equation}
x_1^3+x_2^3+x_3^3-b x_1 x_2 x_3=0,
\end{equation}
where we can make a change of variables:
\begin{align}
\begin{split}
&n_1=\frac{x_2 x_3}{x_1^2}, \quad n_{2}=-\frac{x_1 x_2 x_3 n_3}{x_1 x_2 x_3-(x_2^3+x_3^3) n_3}.
\end{split}
\end{align}

\section{Appendix E}

{\bf Hamiltonian of the Noumi-Yamada system of type $A_4^{(1)}$ (see \cite{{N1,Ta}})}
\begin{align}\label{NY}
\begin{split}
&\frac{dq_i}{dt}=\frac{\partial H_{A_4^{(1)}}}{\partial p_i}, \quad \frac{dp_i}{dt}=-\frac{\partial H_{A_4^{(1)}}}{\partial q_i} \quad (i=1,2),\\
&H_{A_4^{(1)}}(q_1,p_1,q_2,p_2,t;\alpha_0,\alpha_1,\alpha_2,\alpha_3,\alpha_4)\\
&=\frac{t q_1 p_1-q_1^2 p_1-q_1 p_1^2-\alpha_1 q_1+(\alpha_2+\alpha_4) p_1+t q_2 p_2-q_2 p_2^2-q_2^2 p_2-\alpha_3 q_2+\alpha_4 p_2-2p_1 q_2 p_2}{\alpha_0+\alpha_1+\alpha_2+\alpha_3+\alpha_4}.
\end{split}
\end{align}
The system \eqref{NY} admits affine Weyl group symmetry of type $A_4^{(1)}$ as the group of its B{\"a}cklund transformations.

{\bf Holomorphy conditions (cf. \cite{{N1,Ta}})}
\begin{align}\label{holoNY}
\begin{split}
&r_0:x_0=-((q_1+p_1+p_2-t)p_1-\alpha_0)p_1,\ y_0=\frac{1}{p_1}, \ z_0=q_2+p_1, \ w_0=p_2,\\
&r_1:x_1=\frac{1}{q_1}, \ y_1=-(p_1 q_1+\alpha_1)q_1, \ z_1=q_2, \ w_1=p_2,\\
&r_2:x_2=-((q_1-q_2)p_1-\alpha_2)p_1, \ y_2=\frac{1}{p_1}, \ z_2=q_2, \ w_2=p_2+p_1,\\
&r_3:x_3=q_1,\ y_3=p_1, \ z_3=\frac{1}{q_2}, \ w_3=-(p_2 q_2+\alpha_3) q_2,\\
&r_4:x_4=q_1,\ y_4=p_1, \ z_4=-(q_2 p_2-\alpha_4) p_2, \ w_4=\frac{1}{p_2},
\end{split}
\end{align}
where the system \eqref{NY} has the following invariant divisors (cf. \cite{N1,Ta}):
\begin{center}
\begin{tabular}{|c||c|c|} \hline
parameter's relation & $f_i$ \\ \hline
$\alpha_0=0$ & $f_0:=q_1+p_1+p_2-t$  \\ \hline
$\alpha_1=0$ & $f_1:=p_1$  \\ \hline
$\alpha_2=0$ & $f_2:=q_1-q_2$  \\ \hline
$\alpha_3=0$ & $f_3:=p_2$  \\ \hline
$\alpha_4=0$ & $f_4:=q_2$  \\ \hline
\end{tabular}
\end{center}
The B{\"a}cklund transformations of this system satisfy
\begin{equation}
s_i(g)=g+\frac{\alpha_i}{f_i}\{f_i,g\}+\frac{1}{2!} \left(\frac{\alpha_i}{f_i} \right)^2 \{f_i,\{f_i,g\} \}+\cdots \quad (g \in {\mathbb C}(t)[q_1,p_1,q_2,p_2]),
\end{equation}
where $\{,\}$ is the Poisson bracket such that $\{p_i,q_j\}={\delta}_{ij}$ (cf. \cite{N1}).

{\bf Relation of principal eigenvalues $n_i$}
\begin{align}
\begin{split}
&n_1 n_2 n_3 n_4 n_5-n_1 n_2 n_3-n_1 n_2 n_5-n_1 n_4 n_5-n_2 n_3 n_4-n_3 n_4 n_5\\
&+n_1+n_2+n_3+n_4+n_5-2=0,
\end{split}
\end{align}
where this case is $n_i=2 \ (i=1,2,\ldots,5)$. This equation is invariant under ${\mathbb Z}/5{\mathbb Z}$.

\end{document}